\documentclass[final]{article}

\usepackage{graphicx}
\usepackage{amsfonts}
\usepackage{amsmath}
\usepackage{amssymb}
\usepackage{fancyhdr}
\usepackage{titlesec}
\usepackage{indentfirst}
\usepackage{booktabs}
\usepackage{verbatim}
\usepackage{color}
\usepackage{latexsym}
\usepackage{amscd}
\usepackage{esvect}
\usepackage{graphicx}
\usepackage{upgreek}
\usepackage{subcaption}   
\usepackage{caption}
\usepackage[page,header]{appendix}
\usepackage{titletoc}
\usepackage{mathrsfs}
\usepackage[numbers]{natbib}
\usepackage{float}
\usepackage{lipsum}
\usepackage[export]{adjustbox} 
\numberwithin{equation}{section}
\newtheorem{theorem}{Theorem}[section]

\newtheorem{remark}[theorem]{Remark}

\allowdisplaybreaks

\begin{document}
\title{Asymptotic-preserving schemes for two-species binary collisional kinetic system with disparate masses I: 
time discretization and asymptotic analysis
 \footnote{The first and the third author was supported by the funding DOE--Simulation Center for Runaway Electron Avoidance and Mitigation. The second author was supported by NSF grants DMS-1522184 and DMS-1107291: RNMS KI-Net and NSFC grant No. 31571071. 
}}

\date{\today}
\author{Irene M. Gamba \footnote{Department of Mathematics and Oden Institute for Computational Engineering and Sciences, University of Texas at Austin, Austin, TX 78712, USA (gamba@math.utexas.edu).}, 
Shi Jin \footnote{School of Mathematical Sciences, Institute of Natural Sciences, MOE-LSEC and SHL-MAC,
Shanghai Jiao Tong University, Shanghai 200240, China (shijin-m@sjtu.edu.cn).}
and Liu Liu\footnote{Oden Institute for Computational Engineering and Sciences, University of Texas at Austin, Austin, TX 78712, USA (lliu@ices.utexas.edu).}}
\maketitle

\centerline{\it To the memory of David Shen-ou Cai}
\bigskip

\abstract{We develop efficient asymptotic-preserving time discretization schemes to solve the disparate mass kinetic system of a binary gas or plasma
  in the ``relaxation time scale'' relevant to 
  the epochal relaxation phenomenon. {\color{black} Since the resulting model is associated to a parameter given by the square of the mass ratio
  between the light and heavy particles, 
  we develop an asymptotic-preserving scheme in a novel decomposition strategy using the penalization method for multiscale collisional kinetic     equations. }
  Both the Boltzmann and Fokker-Planck-Landau (FPL) binary collision operators will be considered. 
  Other than utilizing several AP strategies for single-species binary kinetic equations, we also
  introduce a novel splitting and a carefully designed explicit-implicit
  approximation, which are guided by the asymptotic analysis of the system.
  We also conduct asymptotic-preserving analysis for the time discretization, for both space homogenous
  and inhomogeneous systems. }

\section{Introduction}

We are interested in the numerical approximation of a disparate mass binary gas or
plasma system, consisting of the mixture of light particles and the heavy ones. Depending on different scalings, such a mixture exhibits various different and interesting asymptotic behavior which poses tremendous numerical
challenges due to both the strongly coupled collisional mechanism, described
by the nonlinear and nonlocal Boltzmann or Fokker-Planck-Landau (FPL)
collision operators, and multiple time and space scales. In the case of plasma, a mixture of electrons and ions,
the equalization of electron and ion temperatures is one of the oldest problems in plasma
physics and was initially considered by Landau \cite{Landau}. 
See \cite{MC-Bobylev, BPS, GHJ, Bird, Johnson2, Johnson1} for more physical description of gas mixtures. By introducing the small scaling parameter, which is the square root of the ratio between the masses of the two kinds of particles, one can obtain various
interesting asymptotic limits by different time scalings of the equations,
see \cite{BPS, Degond, Degond2} for both the Boltzmann and FPL collisions. In particular, under the so-called
``relaxation time scale",
both particle distribution functions are thermalized and the temperatures evolve toward each other via a relaxation equation. This is the 
the epochal relaxation phenomenon first pointed out by Grad \cite{Grad}, and
is the asymptotic regime we are interested in here.
For recent numerical studies of the disparate mass problems, see
\cite{JinLi, Cheng-Gamba}.

One of the main computational challenges for multiscale kinetic equations for binary interactions
is the necessity to resolve the small, microscopic scales numerically which
are often computationally prohibitive. 
In this regard, the Asymptotic-Preserving (AP) schemes \cite{jin1999efficient} have been very
popular in the kinetic and hyperbolic communities in the last two decades.
Such schemes allow one to
use {\it small-scale independent} computational parameters in regimes where one
cannot afford to resolve the small physical scalings numerically.
Such schemes are designed such that they mimic the asymptotic
transition from one scale to another at the discrete level, and
also use specially designed explicit-implicit time discretizations so as
to reduce the algebraic complexity when implicit discretizations are needed.
See review articles \cite{jin2010asymptotic, degond2017asymptotic}. 
For single species particles, in order to overcome the stiffness of the collision
operators, one could penalize the collision operators by simple ones that
are easier to invert, see \cite{Filbet-Jin, JinYan}, or uses exponential
Runge-Kutta methods \cite{dimarco2011exponential, li2014exponential},
or via the micro-macro decomposition \cite{MM-Lemou, Lemou}. See also 
\cite{liu2016unified}.
However, for binary interactions in multispecies models, one encounters extra difficulties due to the
coupling of collision terms between different species. The Cauchy problem for the full non-linear homogeneous Boltzmann system 
describing multi-component monatomic gas mixtures has been studied recently in \cite{MG18}. 
For relatively simpler scalings which lead to hydrodynamic limits, multispecies AP
schemes were developed in for examples \cite{jin2010micro, JinLi, liu2017unified}. See also \cite{MHG14}, where 
a spectral-Lagrangian Boltzmann solver for a multi-energy level gas was developed. 
However, none of the previous works dealt with the disparate mass systems under the long-time scale studied in this paper.

The main challenges to develop efficient AP schemes for the problems under
study include: 1) the strong coupling of the binary collision terms between
different species; 2) the disparate mass scalings so different species
evolve with different time scales thus different species needed to be
treated differently and 3) the long-time scale. In fact, other than utilizing
several existing AP techniques for single species problems, {we also
introduce two {\it new} ideas: {\it a novel splitting} of the system,
guided by the asymptotic analysis introduced in \cite{Degond}, which
 is a natural formulation for the design of AP schemes, 
 and identifying less stiff terms from the stiff ones, again taking
 advantage of the asymptotic behavior of the collision operators.
 We will handle both the Boltzmann and FPL collision terms, thanks to their bilinear structure, 
 and in the end the algebraic complexity, judged by the kind of algebraic
 systems to be inverted, somehow similar to the single species counterparts
 as in \cite{Filbet-Jin} and \cite{JinYan}.

 Due to the complexity of the systems under study, we split our
 results in several papers. In the current paper we focus on the time
 discretization, which is the most difficult part for the design of 
 AP schemes for such a system. We will conduct an AP analysis for a simplified
version of the time
 discretization, as was done for their single-species counterpart in
 \cite{Filbet-Jin}. Given the length of the paper, we will leave
 the numerical experiments in a forthcoming paper.

 This paper is organized as follows. In Section \ref{sec1}, we present the
 physical equations and outline their basic properties and the scalings.
 We also review the asymptotic analysis in \cite{Degond} for the space homogenous case, under the relaxation time scaling. In Section \ref{sec2}, an AP time discretization for the space homogeneous 
equations will be presented, with an asymptotic analysis of its AP property.  
Section \ref{sec3} extends the scheme and analysis to the space inhomogeneous case, by combining with the idea of diffusive relaxation schemes in \cite{JPT1, JPT2} to handle the (also stiff) convection terms.
Conclusions and future works will be given in Section \ref{sec4}. 

\section{An Overview}
\label{sec1}

In this section we present the physical equations which include both Boltzmann and FPL collisions, 
their scalings and fundamental properties, and the asymptotic limit conducted 
in \cite{Degond}. 

\subsection{The equations and scalings}

Let $f^L(t,x,v)$ and $f^H(t,x,v)$ be the probability density distributions of the light and heavy particles at time $t$, position $x$ with 
velocity $v$. 
The rescaled, space inhomogeneous equations are given by
\begin{align}
&\displaystyle\label{scale1} \frac{\partial f^L}{\partial t} + v^L\cdot\nabla_x f^L + F^L \cdot\nabla_{v^L}f^L 
= \mathcal Q^{LL}(f^L, f^L) + \mathcal Q^{LH}_{\varepsilon}(f^L, f^H), \\[2pt]
&\displaystyle\label{scale2} \frac{\partial f^H}{\partial t}  + \varepsilon \left(v^H \cdot\nabla_x f^H + F^H \cdot\nabla_{v^H}f^H\right)
= \varepsilon\left[\mathcal Q^{HH}(f^H, f^H) + \mathcal Q^{HL}_{\varepsilon}(f^H, f^L)\right], 
\end{align}
where $F^L$, $F^H$ stand for the force fields. The definitions of collision operators $\mathcal Q^{LL}$, $\mathcal Q^{HH}$, $\mathcal Q^{LH}_{\varepsilon}$ and $\mathcal Q^{HL}_{\varepsilon}$ represent the binary collisions between light (`L') and heavy
(`H') particles, are given in the Appendix, since only some of their properties, not their specific forms, will be used in this paper. 
Moreover, we assume these are binary interaction operators with transition probability rates presenting the natural symmetries that give rise to 
the classical conservation laws for mixtures. 
$\varepsilon$ is the square root of the mass ratio between the light and heavy particles. 

Define $n$, $u$ and $T$ the density, bulk velocity, and temperature                               
\begin{equation}\label{Rho} n = \int_{\mathbb R^3}\, f(v)\, dv, \qquad u = \frac{1}{n}\, \int_{\mathbb R^3}\, f(v)v\, dv, 
\qquad T = \frac{1}{3 \rho}\, \int_{\mathbb R^3}\, f(v)|v-u|^2\, dv, 
\end{equation}
and denote $M_{u,T}$ the normalized Maxwellian
\begin{equation}\label{normal_M}
M_{u,T}(v)=\frac{1}{(2\pi T)^{3/2}}\,\exp\left(-\frac{|v-u|^2}{2T}\right). 
\end{equation}

In \cite{Degond}, three different time scales were introduced which lead to different hydrodynamic limits. 
We are interested in the third time scale, namely the ``relaxation time scale" studied in \cite{Degond}. The macroscopic limit under this scaling, as well as the design of AP schemes, are the most challenging. The AP schemes that preserve the other two 
asymptotic limits are easy to design by classical AP strategies so will not be discussed 
here. 

The collision time for the light and heavy species are denoted by $t_0^{L}$ and $t_0^{H}$, respectively. We define $t_0=t_0^{L}$ as the basic time scale. Introduce the long time scaling $t_0^{\prime}=t_0/\varepsilon^2$ 
and change of variables $t^{\prime}=\varepsilon^2 t$, $x^{\prime}=\varepsilon x$, $F^{\prime}=F/\varepsilon$, 
at which both distribution functions will be thermalized and temperatures influence each other via a relaxation equation. 
Then the evolution equations are given by 
\begin{align}
&\displaystyle\label{SI_L} \frac{\partial f^L}{\partial t} + \frac{1}{\varepsilon}\left(v^L\cdot\nabla_x f^L + F^L\cdot\nabla_{v^L}f^L\right) = \frac{1}{\varepsilon^2}\left[\mathcal Q^{LL}(f^L, f^L) + \mathcal Q^{LH}_{\varepsilon}(f^L, f^H)\right], \\[4pt]
&\displaystyle\label{SI_H} \frac{\partial f^H}{\partial t}  + \left(v^H\cdot\nabla_x f^H + F^H\cdot\nabla_{v^H}f^H\right) = \frac{1}{\varepsilon}\left[\mathcal Q^{HH}(f^H, f^H) + \mathcal Q^{HL}_{\varepsilon}(f^H, f^L)\right]. 
\end{align}

Inserting the ansatz
$$\mathcal Q^{LH}_{\varepsilon}=\mathcal Q_0^{LH}+\varepsilon \mathcal Q_1^{LH} + O(\varepsilon^2), \qquad
\mathcal Q^{HL}_{\varepsilon}=\mathcal Q_0^{HL}+\varepsilon \mathcal Q_1^{HL} + O(\varepsilon^2) $$
into (\ref{SI_L})--(\ref{SI_H}), one has
\begin{align}
&\displaystyle\notag
\frac{\partial f^L}{\partial t} + \varepsilon^{-1}\left(v^L\cdot\nabla_x f^L + F^L \cdot\nabla_{v^L}f^L\right) \\[2pt]
&\displaystyle\label{scale3}
= \varepsilon^{-2} \left(\mathcal Q^{LL}(f_{\varepsilon}^L, f_{\varepsilon}^L) + \mathcal Q_0^{LH}(f_{\varepsilon}^L, f_{\varepsilon}^H)\right)
+\varepsilon^{-1}\mathcal Q_1^{LH}(f_{\varepsilon}^L, f_{\varepsilon}^H)+\mathcal Q_2^{LH}(f_{\varepsilon}^L, f_{\varepsilon}^H) + 
O(\varepsilon), \\[4pt]
&\displaystyle\notag
\frac{\partial f^H}{\partial t} + v^H \cdot\nabla_x f^H + F^H \cdot\nabla_{v^H}f^H  \\[2pt]
&\displaystyle\label{scale4}  = \varepsilon^{-1}\left(\mathcal Q^{HH}(f_{\varepsilon}^H, f_{\varepsilon}^H) + \mathcal Q_0^{HL}(f_{\varepsilon}^H, f_{\varepsilon}^L)\right)+\mathcal Q_1^{HL}(f_{\varepsilon}^H, f_{\varepsilon}^L)+O(\varepsilon). 
\end{align} 
{\color{black} Clearly the dynamics of (\ref{scale3})--(\ref{scale4}) have stiff terms associated to the electron-ion mass ratio that 
naturally enables the development of asymptotic-preserving schemes. }
\\[4pt]
\indent We first give a summary of the propositions and lemmas on the properties of the collision operators given in \cite{Chapman58, Petit75} 
and summarized in \cite{Degond} that will be useful in our paper. We call ``inter-particle collisions" and ``intra-particle collisions" 
to distinguish binary collisions between different species and like particles in the sequel. 

\begin{theorem}\ \\
\label{property} 
{\bf 1}. For the FPL collision operator, 
\begin{align}
&\displaystyle\label{Q0_F}\mathcal Q_0^{LH}(f^L, f^H)=n^H\, q_0(f^L), 
\qquad q_0(f^L)=\nabla_{v^L}\cdot\left[B(v^L)S(v^L)\nabla_{v^L}f^L(v^L)\right],  \\[2pt]
&\displaystyle\notag\mathcal Q_0^{HL}(f^H, f^L)=-2\,\nabla_{v^H}f^H(v^H)\cdot\int_{\mathbb R^3}\, \frac{B(v^L)}{|v^L|^2}v^L f^L(v^L)\, dv^L. 
\end{align}
For the Boltzmann collision operator, 
\begin{align}
&\displaystyle\label{Q0_B} \mathcal Q_0^{LH}(f^L, f^H) = n^H\, q_0(f^L), \,
q_0(f^L) = \int_{\mathbb S^2} B(v^L, \Omega)\left(f^L(v^L - 2(v^L, \Omega)\Omega) - f^L(v^L)\right)d\Omega, \\[2pt]
&\displaystyle\notag\mathcal Q_0^{HL}(f^H, f^L) = -2 \nabla_{v^H}f^H \cdot 
\int_{\mathbb R^3\times\mathbb S^2} B(v^L, \Omega)\, \frac{(v^L, \Omega)^2}{|v^L|^2}\, v^L f^L(v^L)dv^L d\Omega. 
\end{align}

{\bf 2}. For any function $f^H$, \\ 
\indent (i) if $f^L$ is a function of $|v^L|$, then $\mathcal Q_0^{LH}(f^L, f^H)=0$, \\[2pt]
\indent (ii) if $f^L$ is an even function, then $\mathcal Q_0^{HL}(f^H, f^L)=0$. 

{\bf 3}. The conservation properties of the inter-particle collision operators are given by 
\begin{align}
&\displaystyle\label{Q_e} \int_{\mathbb R^3} \mathcal Q^{LH}_{\varepsilon}\, dv^L = \int_{\mathbb R^3} \mathcal Q^{HL}_{\varepsilon}\, dv^H=0, \\[2pt]
&\displaystyle\notag \int_{\mathbb R^3}\mathcal Q^{LH}_{\varepsilon}\begin{pmatrix}v^L \\ |v^L|^2 \end{pmatrix} dv^L 
+ \varepsilon  \int_{\mathbb R^3} \mathcal Q^{HL}_{\varepsilon}\begin{pmatrix}v^H \\ |v^H|^2 \end{pmatrix} dv^H = 0, \\[2pt]
&\displaystyle\label{Q_i} \int_{\mathbb R^3} \mathcal Q_i^{LH}\, dv^L = \int_{\mathbb R^3} \mathcal Q_i^{HL}\, dv^H=0, \qquad \forall i\in\mathbb N,  \\[2pt]
&\displaystyle\notag \int_{\mathbb R^3} \mathcal Q_i^{LH}\, v^L \, dv^L +
\int_{\mathbb R^3} \mathcal Q_i^{HL}\, v^H\, dv^H=0, \qquad\forall i\in\mathbb N, \\[2pt]
&\displaystyle\notag \int_{\mathbb R^3} \mathcal Q_i^{LH}\, |v^L|^2\, dv^L 
+ \int_{\mathbb R^3} \mathcal Q_{i-1}^{HL}\, |v^H|^2\, dv^H = 0, \qquad\forall i\in\mathbb N, \, i\geq 1, \\[2pt]
&\displaystyle\label{Q0_LH}  \int_{\mathbb R^3} \mathcal Q_0^{LH}\, |v^L|^2\, dv^L =0. 
\end{align}

{\bf 4}. Introduce the operator 
$$\mathcal Q_0^L(f^L) := \mathcal Q^{LL}(f^L, f^L) + n^H q_0(f^L). $$
For all sufficiently regular $f$, $$\int_{\mathbb R^3} \mathcal Q_0^L(f^L)\, \ln f\, dv \leq 0, $$
and
\begin{equation}\label{Q0_L}\mathcal Q_0^L(f)=0  \Leftrightarrow \exists (n, T) \in [0, \infty)^2 \text{ such that }
f = n M_{0,T}, \end{equation}
where $M_{0,T}$ is the normalized Maxwellian defined in (\ref{normal_M}) with $u=0$. 

{\bf 5}. {\color{black}Define $M_0^L:=n_0^L(t)\, M_{0, T_0^L(t)}$.} $\Gamma_0^L$ is a non-positive self-adjoint operator associated with the inner product
$$\langle \phi, \psi \rangle = \int_{\mathbb R^3} \phi\, \psi\, M_0^L\, dv$$
on the space $\chi=\{\phi(v), \langle \phi, \phi \rangle< \infty \}$, and is such that 
$$ \text{ker } \Gamma_0^L = \{ \phi(v^L) \text{ such that } 
\exists (a, b)\in \mathbb R^2, \, \phi(v^L)=a+b|v^L|^2 \}. $$
For $\psi\in \chi$, the equation $\Gamma_0^L \phi = \psi$ 
is solvable if and only if 
$$ \int_{\mathbb R^3} \psi\begin{pmatrix} 1 \\ 
|v^L|^2 \end{pmatrix} M_0^L\, dv^L=0. $$ 
Then the solution $\phi$ is unique in $\left(\text{ker }\Gamma_0^L
\right)^{\perp}$.
\end{theorem}

\subsection{The macroscopic approximation}
\label{Macro}

For clarity of the presentation, we first consider the space homogeneous case of (\ref{scale3})--(\ref{scale4}), so the spatial and velocity gradients on the left-hand-side of the equations are omitted. 
Inserting the Hilbert expansions
$$f_{\varepsilon}^L = f_0^L + \varepsilon f_1^L + \varepsilon^2 f_2^L+\cdots, \qquad
f_{\varepsilon}^H =f_0^H +\varepsilon f_1^H +\varepsilon^2 f_2^H+\cdots $$
and equating terms of $\varepsilon$ leads to: \\[2pt]
 
\noindent order $\varepsilon^{-2}$: 
\begin{equation}\label{e2}\mathcal Q^{LL}(f_0^L, f_0^L)+\mathcal Q_0^{LH}(f_0^L, f_0^H)=0; \end{equation}
order $\varepsilon^{-1}$: 
\begin{align}
&\displaystyle \label{e1_1} 0 = 2\mathcal Q^{LL}(f_0^L, f_1^L) +\mathcal Q_0^{LH}(f_0^L, f_1^H) + \mathcal Q_0^{LH}(f_1^L, f_0^H) +
\mathcal Q_1^{LH}(f_0^L, f_0^H), \\[4pt]
&\displaystyle \label{e1_2} 0 = \mathcal Q^{HH}(f_0^H, f_0^H) +\mathcal Q_0^{HL}(f_0^H, f_0^L); 
\end{align}
order $\varepsilon^0$: 
\begin{align}
&\displaystyle  \frac{\partial f_0^L}{\partial t} = 2\mathcal Q^{LL}(f_0^L, f_2^L)+ \mathcal Q^{LL}(f_1^L, f_1^L) +\mathcal Q_0^{LH}(f_0^L, f_2^H) 
+\mathcal Q_0^{LH}(f_1^L, f_1^H)+\mathcal Q_0^{LH}(f_2^L, f_0^H) \notag\\[4pt]
&\displaystyle \label{e0_1} \qquad\quad +\mathcal Q_1^{LH}(f_0^L, f_1^H) + \mathcal Q_1^{LH}(f_1^L, f_0^H) +
\mathcal Q_2^{LH}(f_0^L, f_0^H),  \\[4pt]
&\displaystyle  \label{e0_2}   \frac{\partial f_0^H}{\partial t} =2\mathcal Q^{HH}(f_0^H, f_1^H) +\mathcal Q_0^{HL}(f_0^H, f_1^L) +\mathcal Q_0^{HL}(f_1^H, f_0^L)
+\mathcal Q_1^{HL}(f_0^H, f_0^L). 
\end{align}
\\[2pt]
First consider the equation for the heavy particles. 
By (\ref{Q0_F}) and (\ref{Q0_B}), we know $$\mathcal Q_0^{LH}(f^L, f^H) = n_0^H\, q_0(f^L), $$ 
with different $q_0(f^L)$ definitions for the Boltzmann and FPL equations respectively. 
Using (\ref{Q0_L}), equation (\ref{e2}) gives 
$f_0^L=M_0^L$. 
By statement 2(ii) in Theorem \ref{property}, since $f_0^L$ is an even function, thus 
$$\mathcal Q_0^{HL}(f_0^H, f_0^L)=0, $$ and (\ref{e1_2}) reduces to 
$$\mathcal Q^{HH}(f_0^H, f_0^H)=0. $$
Using the classical theory of the Boltzmann equation \cite{CC}, 
$\exists\, (n_0^H(t), T_0^H(t))\in [0, \infty)^2$, $u_0^H(t)\in\mathbb R^3$, such that 
$$f_0^H=n_0^H(t)\, M_{u_0^H(t), T_0^H(t)}:=M_0^H. $$

By statement 2(i) in Theorem \ref{property}, 
$$\mathcal Q_0^{LH}(f_0^L, f_1^H)=0, $$ since $f_0^L=M_0^L$ is a function of $|v^L|$. Then
(\ref{e1_1}) is an equation for $f_1^L$, which can be solved by setting $$\phi_1^L=f_1^L\, (M_0^L)^{-1}$$ and 
\begin{equation}\label{Gam}\Gamma_0^L \phi_1^L = -(M_0^L)^{-1}\, \mathcal Q_1^{LH}(M_0^L, f_0^H), \end{equation}
where $\Gamma_0^L$ is an operator defined by
\begin{equation}\label{Gam0}\Gamma_0^L \phi = (M_0^L)^{-1}\left[2\mathcal Q^{LL}(M_0^L, M_0^L\phi)+n_0^{H}\, q_0(M_0^L\phi)\right]. 
\end{equation}
According to statement 5 in Theorem \ref{property}, $\Gamma_0^L\phi=\psi$ is solvable if and only if 
\begin{equation}\label{solvability}\int_{\mathbb R^3}\, \psi \begin{pmatrix} 1 \\ |v^L|^2 \end{pmatrix} M_0^L\, dv^L=0. \end{equation}
Therefore, we have $$\psi = -(M_0^L)^{-1}\, \mathcal Q_1^{LH}(M_0^L, f_0^H)$$ in (\ref{Gam}), and (\ref{solvability}) is satisfied thanks to 
statement 5 in Theorem \ref{property}, thus (\ref{Gam}) is solvable and its unique solution 
in $(\text{ker}\Gamma_0^L)^{\perp}$ is given by
$$f_1^L(v^L)=\frac{1}{T_0^L}\, M_0^L(v^L)u_0^H\cdot v^L. $$ 

Since again $\mathcal Q_0^{HL}(f_1^H, f_0^L)=0$, (\ref{e0_2}) is an equation for $f_1^H$, which can be written in terms of 
$\phi^H = f_1^H\, (M_0^H)^{-1}$ with 
\begin{equation}\label{Gam1}\Gamma_0^H \phi^H =
(M_0^H)^{-1}\left[\frac{\partial M_0^H}{\partial t}-\mathcal Q_0^{HL}(M_0^H, f_1^L)-\mathcal Q_1^{HL}(M_0^H, M_0^L)\right],
 \end{equation}
where $\Gamma_0^H$ is the linearization of $\mathcal Q^{HH}$ around a Maxwellian $M_0^H$:
$$\Gamma_0^H \phi = 2 (M_0^H)^{-1}\, \mathcal Q^{HH}(M_0^H, M_0^H\phi)\,.$$
The necessary and sufficient condition of solvability of equation (\ref{Gam1}) is given by
\begin{equation}\label{solvability2}\int_{\mathbb R^3}\, 
\left[\frac{\partial M_0^H}{\partial t}-\mathcal Q_0^{HL}(M_0^H, f_1^L)-\mathcal Q_1^{HL}(M_0^H, M_0^L)\right] 
\begin{pmatrix}1 \\ v^H \\ |v^H|^2 \end{pmatrix} dv^H=0. \end{equation}
The calculation in \cite{Degond} gives 
$$\int_{\mathbb R^3}\, \left[\mathcal Q_0^{HL}(M_0^H, f_1^L) + \mathcal Q_1^{HL}(M_0^H, M_0^L)\right]
\begin{pmatrix}1 \\ v^H \\ \frac{|v^H|^2}{2} \end{pmatrix} dv^H =
\begin{pmatrix} 0 \\ 0 \\ -3\frac{\lambda(T_0^L)}{T_0^L}\, n_0^L\, n_0^H (T_0^H - T_0^L)\end{pmatrix}\,. $$
Inserting it into (\ref{solvability2}), one finally has
$$\frac{d}{dt}
\begin{pmatrix} n_0^H \\ n_0^H\, u_0^H\\ n_0^H (\frac{1}{2}|u_0^H|^2 +\frac{3}{2}T_0^H)\end{pmatrix}
=\begin{pmatrix} 0 \\ 0 \\  -3\frac{\lambda(T_0^L)}{T_0^L}\, n_0^L\, n_0^H (T_0^H - T_0^L)\end{pmatrix}. $$
Therefore the macroscopic limit of the heavy particles, as $\varepsilon\to 0$, is 
\begin{align*}
&\displaystyle \frac{d}{dt}n_0^H=0, \qquad \frac{d}{dt}n_0^H\, u_0^H=0, \\[2pt]
&\displaystyle \frac{d}{dt}\left(\frac{3n_0^H T_0^H}{2}\right) = -3 \frac{\lambda(T_0^L)}{T_0^L}\, n_0^L n_0^H (T_0^H - T_0^L). 
\end{align*}

Now we consider the light particles. Equation (\ref{e0_2}) is an equation of $f_2^L$ which can be written in terms of $\phi_2^L=f_2^L\, (M_0^L)^{-1}$ with 
$$\Gamma_0^L \phi_2^L =(M_0^L)^{-1}\, S^L, $$
where $\Gamma_0^L$ is defined by (\ref{Gam0}) and
\begin{align}
&\displaystyle S^L = \frac{\partial M_0^L}{\partial t} -\mathcal Q^{LL}(f_1^L, f_1^L) -\mathcal Q_0^{LH}(f_1^L, f_1^H) \notag\\[4pt]
&\displaystyle \label{SL}\qquad -\mathcal Q_1^{LH}(M_0^L, f_1^H)-\mathcal Q_1^{LH}(f_1^L, M_0^H) -\mathcal Q_2^{LH}(M_0^L, M_0^H). 
\end{align}
According to statement 5 in Theorem \ref{property}, the necessary and sufficient condition for the existence of $f_2^L$ should be
\begin{equation}\label{f2_L}\int_{\mathbb R^3}\, S^L(v^L)\begin{pmatrix} 1\\ |v^L|^2\end{pmatrix} dv^L=0\,. \end{equation}
The first equation leads to $dn_0^L/dt=0$. 
By statement 3 in Theorem \ref{property}, 
\begin{align*}
&\displaystyle  \int \mathcal Q_0^{LH}(f_1^L, f_1^H)\, |v^L|^2\, dv^L=0, \\[2pt]
&\displaystyle  \int \mathcal Q_1^{LH}(M_0^L, f_1^H)\, |v^L|^2\, dv^L =\int \mathcal Q_0^{HL}(f_1^H, M_0^L)\, |v^H|^2\, dv^H=0. 
\end{align*}
The remaining terms on the right-hand-side of (\ref{SL}) give 
\begin{align*}
&\displaystyle \int \left[\mathcal Q_1^{LH}(f_1^L, M_0^H) + \mathcal Q_2^{LH}(M_0^L, M_0^H)\right] |v^L|^2\, dv^L \\[2pt]
&\displaystyle \qquad = -\int\left[\mathcal Q_0^{HL}(M_0^H, f_1^L) +\mathcal Q_1^{HL}(M_0^H, M_0^L)\right] |v^H|^2\, dv^H \\[2pt]
&\displaystyle \qquad  = 6 \frac{\lambda(T_0^L)}{T_0^L}\, n_0^L n_0^H\, (T_0^H - T_0^L)\,. 
\end{align*}
Inserting into (\ref{f2_L}), one obtains the evolution equation for $T_0^L$:
\begin{equation}\frac{d}{dt}\left(\frac{3 n_0^L T_0^L}{2}\right) = -3 \frac{\lambda(T_0^L)}{T_0^L}\, n_0^L n_0^H (T_0^L - T_0^H)\,.
\end{equation}
We now summarize the macroscopic equations for the whole system, as $\varepsilon\to 0$, 
\begin{align}
\label{T2}
\begin{split}
&\displaystyle\frac{d}{dt} n_0^L=0,  \\[2pt]
&\displaystyle\frac{d}{dt}\left(\frac{3 n_0^L T_0^L}{2}\right) = -3 \frac{\lambda(T_0^L)}{T_0^L}\, n_0^L n_0^H (T_0^L - T_0^H); \\[2pt]
\end{split}
\end{align}
\begin{align}
\label{T1}
\begin{split}
&\displaystyle\frac{d}{dt} n_0^H = 0, \\[2pt]
&\displaystyle\frac{d}{dt} (n_0^H u_0^H) = 0, \\[2pt]
&\displaystyle\frac{d}{dt}\left(\frac{3n_0^H T_0^H}{2}\right) = -3 \frac{\lambda(T_0^L)}{T_0^L}\, n_0^L n_0^H (T_0^H - T_0^L). 
\end{split}
\end{align}

\section{An asymptotic-preserving time discretization}
\label{sec2}

An AP scheme requires the discrete version of (\ref{SI_L})--(\ref{SI_H}) asymptotically approaches to the macroscopic equations
(\ref{T2})--(\ref{T1}) as $\varepsilon\to 0$, when numerical parameters are held fixed. 
A necessary requirement for such a scheme is some implicit time discretization 
for the numerical stiff terms, which can be easily inverted \cite{jin2010asymptotic}. 
In this section, we design such a time discretization for the space homogeneous equations. 

The space homogeneous version of equations (\ref{SI_L})--(\ref{SI_H}) is given by
\begin{align}
&\displaystyle\label{scale1a}\frac{\partial f^L}{\partial t} = \frac{1}{\varepsilon^2}\left[\mathcal Q^{LL}(f^L, f^L) + \mathcal Q^{LH}_{\varepsilon}(f^L, f^H)\right],  \\[2pt]
&\displaystyle\label{scale2a}\frac{\partial f^H}{\partial t} = \frac{1}{\varepsilon}\left[\mathcal Q^{HH}(f^H, f^H) + {\color{black}\mathcal Q^{HL}_{\varepsilon}}(f^H, f^L)\right]. 
\end{align}

\subsection{A splitting of the equation}

We first decompose $f$ into $f_0$ and $f_1$, 
\begin{equation}f^L=f_0^L+\varepsilon f_1^L, \qquad f^H =f_0^H+\varepsilon f_1^H, 
\end{equation}
and insert into the system (\ref{scale1a})--(\ref{scale2a}), then 
\begin{align}
\label{FL}
\begin{split}
&\displaystyle \frac{\partial}{\partial t}(f_0^L+\varepsilon f_1^L) =\frac{1}{\varepsilon^2}\left[\mathcal Q^{LL}(f_0^L+\varepsilon f_1^L, f_0^L+\varepsilon f_1^L)   
+\mathcal Q^{LH}_{\varepsilon}(f_0^L+\varepsilon f_1^L, f_0^H+\varepsilon f_1^H)\right]\\[2pt]
&\displaystyle \qquad\qquad\qquad = \frac{1}{\varepsilon^2}\bigg[\mathcal Q^{LL}(f_0^L, f_0^L) +2\varepsilon \mathcal Q^{LL}(f_0^L, f_1^L) +\varepsilon^2 \mathcal Q^{LL}(f_1^L, f_1^L) \\[2pt]
&\displaystyle\qquad\qquad\qquad\qquad +\mathcal Q^{LH}_{\varepsilon}(f_0^L, f_0^H) +\varepsilon \mathcal Q^{LH}_{\varepsilon}(f_0^L, f_1^H) +\varepsilon \mathcal Q^{LH}_{\varepsilon}(f_1^L, f_0^H)
+\varepsilon^2 \mathcal Q^{LH}_{\varepsilon}(f_1^L, f_1^H)\bigg], 
\end{split}
\end{align}
and 
\begin{align}
\label{FH}
\begin{split}
&\displaystyle \frac{\partial}{\partial t}(f_0^H+\varepsilon f_1^H) =\frac{1}{\varepsilon}
\left[\mathcal Q^{HH}(f_0^H+\varepsilon f_1^H, f_0^H+\varepsilon f_1^H)   
+ {\color{black}\mathcal Q^{HL}_{\varepsilon}}(f_0^H+\varepsilon f_1^H, f_0^L+\varepsilon f_1^L)\right] \\[2pt]
&\displaystyle \qquad\qquad\qquad = \frac{1}{\varepsilon}\bigg[\mathcal Q^{HH}(f_0^H, f_0^H) +2\varepsilon \mathcal Q^{HH}(f_0^H, f_1^H) + \varepsilon^2 \mathcal Q^{HH}(f_1^H, f_1^H) \\[2pt]
&\displaystyle\qquad\qquad\qquad\qquad + {\color{black}\mathcal Q^{HL}_{\varepsilon}}(f_0^H, f_0^L) +\varepsilon \mathcal Q^{HL}_{\varepsilon}(f_0^H, f_1^L) +\varepsilon \mathcal Q^{HL}_{\varepsilon}(f_1^H, f_0^L)
+\varepsilon^2 \mathcal Q^{HL}_{\varepsilon}(f_1^H, f_1^L)\bigg]. 
\end{split}
\end{align}

Our first key idea is to split (\ref{FL}) into two equations for $f_0^L$, $f_1^L$ respectively, 
\begin{align}
\label{FL-a}
\begin{split}
&\displaystyle\frac{\partial}{\partial t}f_0^L =\frac{1}{\varepsilon^2}\bigg[\mathcal Q^{LL}(f_0^L, f_0^L)+\mathcal Q_0^{LH}(f_0^L, f_0^H)\bigg], \\[4pt]
&\displaystyle \frac{\partial}{\partial t}f_1^L = \frac{1}{\varepsilon^2}\bigg[\frac{1}{\varepsilon}\left(\mathcal Q^{LH}_{\varepsilon}(f_0^L, f_0^H)-\mathcal Q_0^{LH}(f_0^L, f_0^H)\right) +2\mathcal Q^{LL}(f_0^L, f_1^L)+\varepsilon\mathcal Q^{LL}(f_1^L, f_1^L) \\[2pt]
&\displaystyle\qquad\qquad +\mathcal Q^{LH}_{\varepsilon}(f_0^L, f_1^H) +
\mathcal Q^{LH}_{\varepsilon}(f_1^L, f_0^H) +\varepsilon\mathcal Q^{LH}_{\varepsilon}(f_1^L, f_1^H)\bigg], 
\end{split}
\end{align}
and split (\ref{FH}) into two equations for $f_0^H$, $f_1^H$ respectively: 
\begin{align} 
\label{FH-b} 
\begin{split}
&\displaystyle\frac{\partial}{\partial t}f_0^H=\frac{1}{\varepsilon}\bigg[\mathcal Q^{H}(f_0^H, f_0^H)+\mathcal Q_0^{HL}(f_0^H, f_0^L)\bigg], 
\\[4pt]
&\displaystyle\frac{\partial}{\partial t}f_1^H = \frac{1}{\varepsilon}\bigg[\frac{1}{\varepsilon}\left(\mathcal Q^{HL}_{\varepsilon}(f_0^H, f_0^L)-\mathcal Q_0^{HL}(f_0^H, f_0^L)   \right) +2\mathcal Q^{HH}(f_0^H, f_1^H)+\varepsilon \mathcal Q^{HH}(f_1^H, f_1^H) \\[2pt]
&\displaystyle\qquad\qquad +\mathcal Q^{HL}_{\varepsilon}(f_0^H, f_1^L) +
\mathcal Q^{HL}_{\varepsilon}(f_1^H, f_0^L) + \varepsilon \mathcal Q^{HL}_{\varepsilon}(f_1^H, f_1^L)\bigg]. 
\end{split}
\end{align}
This splitting is motivated by the asymptotic analysis presented in subsection \ref{Macro}, and plays the central role in the AP 
time discretization, which will be introduced in the next subsection. 

\subsection{Time discretization}
\label{Time-D}
First, to have a scheme uniformly stable with respect to $\varepsilon$, it is natural to use the implicit discretizations for all the stiff collision terms, namely, those that appear to be of $O(1)$ inside the brackets on the right hand side of (\ref{FL-a})--(\ref{FH-b}). 
We use the notations $f_{L,0}^n$, $f_{L,1}^n$, $f_{H,0}^n$, $f_{H,1}^n$ to denote the numerical solutions of 
$f_0^L$, $f_1^L$, $f_0^H$ and $f_1^H$ at time step $t^n$. 
Consider the light particles. A naive discretization for $f_{L,0}$, $f_{L,1}$ in (\ref{FL-a}) is 
\begin{align}
&\displaystyle \label{num1}\frac{f_{L,0}^{n+1}-f_{L,0}^n}{\Delta t} =\frac{1}{\varepsilon^2}\bigg[\mathcal Q^{LL}(f_{L,0}^{n+1}, f_{L,0}^{n+1})+
\mathcal Q_0^{LH}(f_{L,0}^{n+1}, f_{H,0}^{n+1})\bigg], \\[4pt]
&\displaystyle\frac{f_{L,1}^{n+1}-f_{L,1}^n}{\Delta t} = \frac{1}{\varepsilon^2}\bigg[\frac{1}{\varepsilon}\left(\mathcal Q^{LH}_{\varepsilon}(f_{L,0}^{n+1}, f_{H,0}^{n+1})
-\mathcal Q_0^{LH}(f_{L,0}^{n+1}, f_{H,0}^{n+1})\right) \notag\\[4pt]
&\displaystyle\qquad\qquad\qquad\qquad +2 \mathcal Q^{LL}(f_{L,0}^{n+1}, f_{L,1}^{n+1}) +
\varepsilon\mathcal Q^{LL}(f_{L,1}^{n}, f_{L,1}^{n})\notag\\[4pt]
&\displaystyle \label{num2} \qquad\qquad\qquad\qquad + \mathcal Q^{LH}_{\varepsilon}(f_{L,0}^{n+1}, f_{H,1}^{n+1}) +
\mathcal Q^{LH}_{\varepsilon}(f_{L,1}^{n+1}, f_{H,0}^{n+1}) + \varepsilon\mathcal Q^{LH}_{\varepsilon}(f_{L,1}^{n}, f_{H,1}^n)\bigg].  
\end{align}
Consider the time evolution for $f_{H,0}$, $f_{H,1}$. A naive implicit scheme for (\ref{FH-b}) would be: 
\begin{align}
&\displaystyle\label{num3}\frac{f_{H,0}^{n+1}-f_{H,0}^n}{\Delta t} =\frac{1}{\varepsilon}\bigg[\mathcal Q^{HH}(f_{H,0}^{n+1}, f_{H,0}^{n+1})+
\mathcal Q_0^{HL}(f_{H,0}^{n+1}, f_{L,0}^{n+1})\bigg], \\[6pt]
&\displaystyle\frac{f_{H,1}^{n+1}-f_{H,1}^n}{\Delta t} = \frac{1}{\varepsilon}\bigg[\frac{1}{\varepsilon}\left(\mathcal Q^{HL}_{\varepsilon}(f_{H,0}^{n+1}, f_{L,0}^{n+1})
-\mathcal Q_0^{HL}(f_{H,0}^{n+1}, f_{L,0}^{n+1})\right) \notag\\[4pt]
&\displaystyle\qquad\qquad\qquad\qquad
+2 \mathcal Q^{HH}(f_{H,0}^{n+1}, f_{H,1}^{n+1}) + \varepsilon\mathcal Q^{HH}(f_{H,1}^{n}, f_{H,1}^{n})\notag \\[4pt]
&\displaystyle\label{num4} \qquad\qquad\qquad\qquad + \mathcal Q^{HL}_{\varepsilon}(f_{H,0}^{n+1}, f_{L,1}^{n+1}) +
\mathcal Q^{HL}_{\varepsilon}(f_{H,1}^{n+1}, f_{L,0}^{n+1}) + \varepsilon\mathcal Q^{HL}_{\varepsilon}(f_{H,1}^{n}, f_{L,1}^n)\bigg], 
\end{align}
in which the right-hand-side is fully implicit, except the terms that are relatively less stiff due to an extra factor of $\varepsilon$. 
Inverting the above system is algebraically complex due to the nonlinearity, nonlocal nature of the collision operators and the coupling between the two types of particles. 
Our next key idea is to use the asymptotic behavior of the operators to identify those terms that are not stiff. 

\subsubsection{Identifying the less stiff terms}

\indent First, as $\varepsilon\to 0$, 
\begin{equation}\label{fL_0} f_{L,0}^{n+1}\to n_0^L\, M_{0,T_0^L}. \end{equation}
{\color{black}Since $M_{0,T_0^L}$ is a function of $|v^L|$, according to 2(i) in Theorem \ref{property}, 
$$\mathcal Q_0^{LH}(n_0^L\, M_{0,T_0^L}, f_{H,0}^{n+1})=0, $$
thus $$\mathcal Q_0^{LH}(f_{L,0}^{n+1}, f_{H,0}^{n+1})=O(\varepsilon), $$} which is less stiff and can be implemented explicitly. 

Secondly, as $\varepsilon\to 0$, similarly
$$\mathcal Q^{LH}_{\varepsilon}(f_{L,0}^{n+1}, f_{H,1}^{n+1})\to \mathcal Q_{0}^{LH}(f_{L,0}^{n+1}, f_{H,1}^{n+1})=O(\varepsilon),$$ so the corresponding term is less stiff and can also be discretized explicitly. 

For the less stiff terms $Q_0^{LH}(f_{L,0}, f_{H,0})$ and $\mathcal Q^{LH}_{\varepsilon}(f_{L,0}, f_{H,1})$ we treat them explicitly, 
thus our time discretizations for $f_{L,0}$, $f_{L,1}$ are given by 
\begin{align}
&\displaystyle\label{num1N}\frac{f_{L,0}^{n+1}-f_{L,0}^n}{\Delta t} =\frac{1}{\varepsilon^2}\bigg[\mathcal Q^{LL}(f_{L,0}^{n+1}, f_{L,0}^{n+1})+\mathcal Q_0^{LH}(f_{L,0}^{n}, f_{H,0}^{n})\bigg], \\[6pt]
&\displaystyle\frac{f_{L,1}^{n+1}-f_{L,1}^n}{\Delta t} = \frac{1}{\varepsilon^2}\bigg[\frac{1}{\varepsilon}\left(\mathcal Q^{LH}_{\varepsilon}(f_{L,0}^{n+1}, f_{H,0}^{n+1})
-\mathcal Q_0^{LH}(f_{L,0}^{n+1}, f_{H,0}^{n+1})\right) \notag\\[4pt]
&\displaystyle\qquad\qquad\qquad\qquad + 2\mathcal Q^{LL}(f_{L,0}^{n+1}, f_{L,1}^{n+1}) +\varepsilon\mathcal Q^{LL}(f_{L,1}^{n}, f_{L,1}^{n}) \notag\\[4pt]
&\displaystyle\label{num2N}\qquad\qquad\qquad\qquad + \mathcal Q^{LH}_{\varepsilon}(f_{L,0}^{n}, f_{H,1}^{n}) +
\mathcal Q^{LH}_{\varepsilon}(f_{L,1}^{n+1}, f_{H,0}^{n+1}) + \varepsilon\mathcal Q^{LH}_{\varepsilon}(f_{L,1}^{n}, f_{H,1}^n)\bigg]. 
\end{align}
\\[2pt]

Similarly for $f_{H,0}$, $f_{H,1}$, 
we introduce the following time discretizations for $f_{H,0}, f_{H,1}$ by taking advantages of some terms that are actually not stiff: 
\begin{align}
&\displaystyle\label{num3N}\frac{f_{H,0}^{n+1}-f_{H,0}^n}{\Delta t} = \frac{1}{\varepsilon}\bigg[\mathcal Q^{HH}(f_{H,0}^{n+1}, f_{H,0}^{n+1})+
\mathcal Q_0^{HL}(f_{H,0}^{n}, f_{L,0}^{n})\bigg], \\[6pt]
&\displaystyle\frac{f_{H,1}^{n+1}-f_{H,1}^n}{\Delta t} = \frac{1}{\varepsilon}\bigg[\frac{1}{\varepsilon}\left(\mathcal Q^{HL}_{\varepsilon}(f_{H,0}^{n+1}, f_{L,0}^{n+1})-\mathcal Q_0^{HL}(f_{H,0}^{n+1}, f_{L,0}^{n+1})\right)\notag\\[4pt]
&\displaystyle\qquad\qquad\qquad\qquad
+2 {\mathcal Q^{HH}(f_{H,0}^{n+1}, f_{H,1}^{n+1})} +\varepsilon\mathcal Q^{HH}(f_{H,1}^{n}, f_{H,1}^{n})\notag \\[4pt]
&\displaystyle\label{num4N} \qquad\qquad\qquad\qquad + 
\mathcal Q^{HL}_{\varepsilon}(f_{H,0}^{n+1}, f_{L,1}^{n+1}) +
\mathcal Q^{HL}_{\varepsilon}(f_{H,1}^{n}, f_{L,0}^{n}) +\varepsilon\mathcal Q^{HL}_{\varepsilon}(f_{H,1}^{n}, f_{L,1}^n)\bigg], 
\end{align}
where the argument 2(ii) of Theorem \ref{property} is used, that is, 
since $f_{L,0}^{n+1}$ is asymptotically an even function due to (\ref{fL_0}), one has
$$\mathcal Q_0^{HL}(f_{H,0}^{n+1}, f_{L,0}^{n+1})=O(\varepsilon), $$
thus the second-term on the right-hand-side of (\ref{num3}) is not stiff. 
In addition, as $\varepsilon\to 0$, 
$$\mathcal Q_{\varepsilon}^{HL}(f_{H,1}^{n+1}, f_{L,0}^{n+1})\to\mathcal Q_0^{HL}(f_{H,1}^{n+1}, f_{L,0}^{n+1})=O(\varepsilon). $$
Thus the term $\mathcal Q^{HL}_{\varepsilon}(f_{H,1}^{n+1}, f_{L,0}^{n+1})$ in (\ref{num4}) is less stiff and can be approximated explicitly. 

\subsubsection{Handling of the stiff terms}

First, we point out the terms
$\mathcal Q_{\varepsilon}^{HL}(f_{H,0}^{n+1}, f_{L,0}^{n+1})$, $\mathcal Q_0^{HL}(f_{H,0}^{n+1}, f_{L,0}^{n+1})$ and 
$\mathcal Q_{\varepsilon}^{HL}(f_{H,0}^{n+1}, f_{L,1}^{n+1})$ in (\ref{num4N}), although implicit, can be obtained 
{\it explicitly} since $f_{L,0}^{n+1}$, $f_{H,0}^{n+1}$ and $f_{L,1}^{n+1}$ are already computed from (\ref{num1N}), (\ref{num2N})
and (\ref{num3N}). 

\indent Now we take care of the truly stiff and implicit collision terms in schemes (\ref{num1N})--(\ref{num2N}) and (\ref{num3N})--(\ref{num4N}). 
They will be penalized by an operator that can either be inverted analytically (for the case of the Boltzmann collision \cite{Filbet-Jin}) 
or by a Poisson-type solver (for the case of FPL collision \cite{JinYan}). 
\\[8pt]
\noindent {\bf (i)}
For the stiff and nonlinear term 
$\mathcal Q^{LH}_{\varepsilon}(f_{L,1}^{n+1}, f_{H,0}^{n+1})$ in (\ref{num2N}), motivated by \cite{Filbet-Jin, JinYan}, we use 
$\mathcal Q_0^{LH}(f_{L,1}, f_{H,0})$ which is the leading order asymptotically for $\varepsilon$ small, 
as the penalty operator. The rationale for this is that 
$\mathcal Q_0^{LH}(f_{L,1}, f_{H,0})$ is much easier to be inverted than 
$\mathcal Q^{LH}_{\varepsilon}(f_{L,1}, f_{H,0})$, as will be shown below. 
We substitute $\mathcal Q^{LH}_{\varepsilon}(f_{L,1}^{n+1}, f_{H,0}^{n+1})$ in (\ref{num2N}) by 
$$\underbrace{\mathcal Q^{LH}_{\varepsilon}(f_{L,1}^{n}, f_{H,0}^{n}) - \mathcal Q_0^{LH}(f_{L,1}^{n}, f_{H,0}^{n})}_{\text{less stiff}}+\underbrace{\mathcal Q_0^{LH}(f_{L,1}^{n+1}, f_{H,0}^{n+1})}_{\text{stiff}}. $$

Integrate both sides of (\ref{num3N}) in $v^H$, we get $n_0^H$ does not change from $t^n$ to $t^{n+1}$, 
so we will drop its dependence on $n$. Thus
$$\mathcal Q_0^{LH}(f_{L,1}^{n+1}, f_{H,0}^{n+1}) = n_0^H\, q_0(f_{L,1}^{n+1}), $$ with $q_0$ defined in (\ref{Q0_F}) and (\ref{Q0_B}) for the Boltzmann and FPL equations respectively. 
For the FPL case, 
\begin{equation}\label{inver}\mathcal Q_0^{LH}(f_{L,1}^{n+1}, f_{H,0}^{n+1}) = n_0^H\, \nabla_{v^L}\cdot \left[B(v^L) S(v^L)\nabla_{v^L}f_{L,1}^{n+1}(v^L)\right], 
\end{equation}
thus one only needs to invert a linear FP operator. See \cite{JinYan}. 
For the Boltzmann case, 
\begin{align*}
&\displaystyle\mathcal Q_0^{LH}(f_{L,1}^{n+1}, f_{H,0}^{n+1}) = 
n_0^H \int_{\mathbb S^2} B(v^L, \Omega)\left(f_{L,1}^{n+1}(v^L - 2(v^L, \Omega)\Omega) - f_{L,1}^{n+1}(v^L)\right) d\Omega \\[2pt]
&\displaystyle = n_0^H \int_{\mathbb S^2} B(v^L, \Omega)\, f_{L,1}^{n+1}(v^L - 2(v^L, \Omega)\Omega)\, d\Omega - 
n_0^H f_{L,1}^{n+1}(v^L) \int_{\mathbb S^2} B(v^L, \Omega)\, d\Omega, 
\end{align*}
which is still a nonlocal operator. 
We use the linear penalty method \cite{Lemou} to remove the stiffness here, that is, substitute the above term by 
$$  n_0^H \int_{\mathbb S^2} B(v^L, \Omega)\left(f_{L,1}^n(v^L - 2(v^L, \Omega)\Omega
- f_{L,1}^n(v^L)\right) d\Omega - n_0^H \mu f_{L,1}^n(v^L) + n_0^H \mu f_{L,1}^{n+1}(v^L), $$
where $$\mu=\max_{v^L}\int_{\mathbb S^2} B(v^L, \Omega)\, d\Omega. $$ 
(See discussions in Remark \ref{RK} for the use of linear penalty here instead of
the BGK penalty of Filbet-Jin \cite{Filbet-Jin}.)
\\[2pt]

\noindent {\bf (ii)} To deal with the stiff terms $\mathcal Q^{LL}(f_{L,0}^{n+1}, f_{L,0}^{n+1})$ and 
$\mathcal Q^{HH}(f_{H,0}^{n+1}, f_{H,0}^{n+1})$ in (\ref{num1N}) and (\ref{num3N}) respectively, 
the {\color{black}BGK} penalty is used for the Boltzmann collision operators \cite{Filbet-Jin}, 
while a linear Fokker-Planck operator will be used to penalize for the FPL collision case, as done in \cite{JinYan}. 
Take the term $\mathcal Q^{LL}(f_{L,0}^{n+1}, f_{L,0}^{n+1})/\varepsilon^2$ and the Boltzmann equation as an example. 
The idea is to split it into the summation of a stiff, dissipative part and a non-(or less) stiff, non-dissipative part: 
$$ \frac{\mathcal Q^{LL}(f_{L,0}^{n+1}, f_{L,0}^{n+1})}{\varepsilon^2} = \underbrace{\frac{\mathcal Q^{LL}(f_{L,0}^n, f_{L,0}^n) - \mathcal P(f_{L,0}^n)}{\varepsilon^2}}_{\text{less stiff}}  + 
\underbrace{\frac{\mathcal P(f_{L,0}^{n+1})}{\varepsilon^2}}_{\text{stiff}}, $$
with $\mathcal P(f_{L,0})$ a well-balanced relaxation approximation of 
$\mathcal Q^{LL}(f_{L,0}, f_{L,0})$ and defined by
$$\mathcal P(f_{L,0}) := \beta_1 (M_{\{n, u, T\}} - f_{L,0}), 
\qquad \beta_1 = \sup_{v}\left|\frac{\mathcal Q^{LL}(f_{L,0},f_{L,0})}{f_{L,0}-M_{\{n,u,T\}}}\right|, $$
and the {\it local} Maxwellian distribution function is 
\begin{equation}\label{local}
M_{\{n, u, T\}} =\frac{n}{(2\pi T)^{3/2}}\,\exp\left(-\frac{|v-u|^2}{2T}\right), \end{equation}
and $n$, $u$, $T$ are defined in (\ref{Rho}) with $f=f_{L,0}$. How to obtain $n$, $u$, $T$ from the moment systems of 
$f_{L,0}$ and $f_{H,0}$ will be discussed below. 
See the Appendix for more details of the penalization for both the Boltzmann and FPL cases. 
\\[2pt]

\noindent {\bf (iii)} To deal with the nonlinear collision operators $\mathcal Q^{LL}(f_{L,0}^{n+1}, f_{L,1}^{n+1})$ in (\ref{num2}), 
since $f_{L,0}^{n+1}$ is already computed from (\ref{num1N}), this is essentially a linear operator and we use the classical formula \cite{CC}
\begin{equation}\label{Q-LL} \mathcal Q^{LL}(f_{L,0}^{n+1}, f_{L,1}^{n+1}) = \frac{1}{4}\left[\mathcal Q^{LL}(f_{L,0}^{n+1}+f_{L,1}^{n+1}, f_{L,0}^{n+1}+f_{L,1}^{n+1})
- \mathcal Q^{LL}(f_{L,0}^{n+1}-f_{L,1}^{n+1}, f_{L,0}^{n+1}-f_{L,1}^{n+1})\right]. 
\end{equation}
For each collision term on the right-hand-side of (\ref{Q-LL}) that has the same argument, 
we adopt the linear penalty method as mentioned in \cite{Lemou} to serve the purpose of removing the stiffness. 
The reason why the BGK-type penalty method of Filbet-Jin does not work well here will be explained in Remark \ref{RK} below. 
The strategy is to substitute 
$\mathcal Q^{LL}(f_{L,0}^{n+1}, f_{L,1}^{n+1})$ by 
\begin{align}
&\displaystyle\notag \frac{1}{4}\bigg[\mathcal Q^{LL}(f_{L,0}^{n}+f_{L,1}^{n}, f_{L,0}^{n}+f_{L,1}^{n}) + \mu (f_{L,0}^n+f_{L,1}^n)
 - \mu (f_{L,0}^{n+1}+f_{L,1}^{n+1}) \\[2pt]
&\displaystyle\label{LP} - \left(\mathcal Q^{LL}(f_{L,0}^{n}-f_{L,1}^{n}, f_{L,0}^{n}-f_{L,1}^{n}) + \mu (f_{L,0}^n-f_{L,1}^n)
- \mu (f_{L,0}^{n+1}-f_{L,1}^{n+1})\right)\bigg], 
\end{align}
where $\mu$ is chosen sufficiently large. 
For the FPL equation, let 
$$\mu>\frac{1}{2}\max_{v} \lambda(D(g)), $$
where $g=f_{L,0}\pm f_{L,1}$ and $\lambda(D(g))$ is the spectral radius of $D$ defined by 
$$ D(g) = \int_{\mathbb R^3} B^L(v^L - v_1^L)S(v^L-v_1^L)g_1^L\, dv_1^L. $$
For the Boltzmann equation, let $\mu > \mathcal Q^{LL,-}$, 
where we split the operator $\mathcal Q^{LL}$ in (\ref{LP}) as 
 $$\mathcal Q^{LL}(g) =\mathcal Q^{LL,+}(g) - g\, \mathcal Q^{LL,-}(g), $$
with the definitions $g=f_{L,0}\pm f_{L,1}$ and 
 $$\mathcal Q^{LL,+}(g)=\int_{\mathbb R^3}\int_{\mathbb S^2} B^L(v^L - v_1^L, \Omega) g^{\prime, L}g_1^{\prime, L}\, d\Omega dv_1^L, 
\quad \mathcal Q^{LL,-}(g)=\int_{\mathbb R^3}\int_{\mathbb S^2} B^L(v^L - v_1^L, \Omega)g_1^L\, d\Omega dv_1^L. $$
The collision term $\mathcal Q^{HH}(f_{H,0}^{n+1}, f_{H,1}^{n+1})$ in (\ref{num4N}) is dealt in a similar way. 
\\[2pt]

Now with the penalties plugged into (\ref{num1N})--(\ref{num2N}) and (\ref{num3N})--(\ref{num4N}), our scheme becomes
\begin{align}
&\displaystyle\label{num1N-P0}\frac{f_{L,0}^{n+1}-f_{L,0}^n}{\Delta t} =
\frac{1}{\varepsilon^2}\bigg[\mathcal Q^{LL}(f_{L,0}^{n}, f_{L,0}^{n}) - \mathcal P(f_{L,0}^n) + \mathcal P(f_{L,0}^{n+1})
+\mathcal Q_0^{LH}(f_{L,0}^{n}, f_{H,0}^{n})\bigg], \\[4pt]
&\displaystyle\frac{f_{L,1}^{n+1}-f_{L,1}^n}{\Delta t} = \frac{1}{\varepsilon^2}\bigg[\frac{1}{\varepsilon}\left(\mathcal Q^{LH}_{\varepsilon}(f_{L,0}^{n+1}, f_{H,0}^{n+1}) -\mathcal Q_0^{LH}(f_{L,0}^{n+1}, f_{H,0}^{n+1})\right) \notag\\[2pt]
&\displaystyle\qquad\qquad\qquad\quad + \frac{1}{2}\bigg[\mathcal Q^{LL}(f_{L,0}^{n}+f_{L,1}^{n}, f_{L,0}^{n}+f_{L,1}^{n}) + \mu (f_{L,0}^n+f_{L,1}^n) - \mu (f_{L,0}^{n+1}+f_{L,1}^{n+1}) \notag\\[2pt]
&\displaystyle \qquad\qquad\qquad\qquad - \left(\mathcal Q^{LL}(f_{L,0}^{n}-f_{L,1}^{n}, f_{L,0}^{n}-f_{L,1}^{n}) + \mu (f_{L,0}^n-f_{L,1}^n)
- \mu (f_{L,0}^{n+1}-f_{L,1}^{n+1})\right)\bigg] \notag\\[2pt]
&\displaystyle\qquad\qquad\qquad\quad  + \varepsilon\mathcal Q^{LL}(f_{L,1}^{n}, f_{L,1}^{n})  + \mathcal Q^{LH}_{\varepsilon}(f_{L,0}^{n}, f_{H,1}^{n}) + \left(\mathcal Q^{LH}_{\varepsilon}(f_{L,1}^{n}, f_{H,0}^{n}) - \mathcal Q_0^{LH}(f_{L,1}^{n}, f_{H,0}^{n})\right) \notag\\[2pt]
&\displaystyle\label{num2N-P0}\qquad\qquad\qquad\quad + \mathcal Q_0^{LH}(f_{L,1}^{n+1}, f_{H,0}^{n+1}) +\varepsilon\mathcal Q^{LH}_{\varepsilon}(f_{L,1}^{n}, f_{H,1}^n)\bigg]; 
\end{align}
\begin{align}
&\displaystyle\label{num3N-P0}\frac{f_{H,0}^{n+1}-f_{H,0}^n}{\Delta t} = 
\frac{1}{\varepsilon}\bigg[\mathcal Q^{HH}(f_{H,0}^{n+1}, f_{H,0}^{n+1}) - \mathcal P(f_{H,0}^n) + \mathcal P(f_{H,0}^{n+1})
+ \mathcal Q_0^{HL}(f_{H,0}^{n}, f_{L,0}^{n})\bigg], \\[4pt]
&\displaystyle\frac{f_{H,1}^{n+1}-f_{H,1}^n}{\Delta t} = \frac{1}{\varepsilon}\bigg[\frac{1}{\varepsilon}\left(\mathcal Q^{HL}_{\varepsilon}(f_{H,0}^{n+1}, f_{L,0}^{n+1})-\mathcal Q_0^{HL}(f_{H,0}^{n+1}, f_{L,0}^{n+1})\right)\notag\\[2pt]
&\displaystyle\qquad\qquad\qquad\quad
+ \frac{1}{2}\bigg[\mathcal Q^{HH}(f_{H,0}^{n}+f_{H,1}^{n}, f_{H,0}^{n}+f_{H,1}^{n}) + \mu(f_{H,0}^n+f_{H,1}^n)  
- \mu(f_{H,0}^{n+1}+f_{H,1}^{n+1}) \notag\\[2pt]
&\displaystyle\qquad\qquad\qquad\qquad  - \left(\mathcal Q^{HH}(f_{H,0}^{n}-f_{H,1}^{n}, f_{H,0}^{n}-f_{H,1}^{n}) + \mu(f_{H,0}^n-f_{H,1}^n) 
- \mu(f_{H,0}^{n+1}-f_{H,1}^{n+1})\right)\bigg] \notag \\[2pt]
&\displaystyle\label{num4N-P0} \qquad\qquad\qquad\quad  + \varepsilon\mathcal Q^{HH}(f_{H,1}^{n}, f_{H,1}^{n}) + \mathcal Q^{HL}_{\varepsilon}(f_{H,0}^{n+1}, f_{L,1}^{n+1}) + \mathcal Q^{HL}_{\varepsilon}(f_{H,1}^{n}, f_{L,0}^{n}) + \varepsilon\mathcal Q^{HL}_{\varepsilon}(f_{H,1}^{n}, f_{L,1}^n)\bigg]. 
\end{align}


\begin{remark}
\label{RK}
In this remark, we will explain why the BGK- or Fokker-Planck type penalties do not work well so the linear penalties are used in 
(\ref{num2N-P0}) and (\ref{num4N-P0}). 
One needs to compute the moment systems in order to define the local Maxwellian $M_{\{n,u,T\}}$ in the penalty operators. 
Define the vectors $$ \phi(v^L)=(1, v^L, \frac{|v^L|^2}{2}), \qquad \phi(v^H)=(1, v^H, \frac{|v^H|^2}{2}), $$
and denote \begin{equation}\label{Phi} \phi_1^L = v^L, \qquad \phi_2^L = \frac{|v^L|^2}{2}, \qquad \phi_1^H = v^H, \qquad\phi_2^H = \frac{|v^H|^2}{2}. 
\end{equation}
Denote the moments by 
\begin{equation}\label{Moments}
 n = \int_{\mathbb R^3} f(v)dv := P_0, \qquad n u = \int_{\mathbb R^3} v f(v) dv := P_1, \qquad \int_{\mathbb R^3} \frac{1}{2}|v|^2\, f(v)dv := P_2. 
\end{equation}
Multiplying (\ref{num1N-P0})--(\ref{num2N-P0}) by $\phi(v^L)$, we obtain the moment systems for $f_{L,0}$, $f_{L,1}$:
\begin{align}
&\displaystyle (P_0)_{L,0}^{n+1} = (P_0)_{L,0}^n,  \notag\\[2pt]
&\displaystyle (P_1)_{L,0}^{n+1} = (P_1)_{L,0}^n + \frac{\Delta t}{\varepsilon^2}\int_{\mathbb R^3} v^L\, 
 \mathcal Q_0^{LH}(f_{L,0}^n, f_{H,0}^n)(v^L)\, dv^L,   \notag\\[2pt]
&\displaystyle (P_2)_{L,0}^{n+1} = (P_2)_{L,0}^n,  \notag\\[4pt]
&\displaystyle\label{M_L1_a} (P_0)_{L,1}^{n+1} = (P_0)_{L,1}^n + \frac{\mu\Delta t}{\varepsilon^2}\Big((P_0)_{L,1}^n - (P_0)_{L,1}^{n+1}\Big), \\[2pt]
&\displaystyle (P_1)_{L,1}^{n+1} = (P_1)_{L,1}^n + \frac{\Delta t}{\varepsilon^2}\int_{\mathbb R^3}\bigg[\frac{1}{\varepsilon}
\left(\mathcal Q_{\varepsilon}^{LH}(f_{L,0}^{n+1}, f_{H,0}^{n+1})(v^L) - \mathcal Q_0^{LH}(f_{L,0}^{n+1}, f_{H,0}^{n+1})(v^L)\right) \notag\\[2pt]
&\displaystyle \qquad\qquad\qquad\qquad + \left(\mathcal Q_{\varepsilon}^{LH}(f_{L,0}^n, f_{H,1}^n)(v^L) + \mathcal Q_{\varepsilon}^{LH}(f_{L,1}^n, f_{H,0}^n)(v^L) - \mathcal Q_0^{LH}(f_{L,1}^n, f_{H,0}^n)(v^L)\right) \notag\\[2pt]
&\displaystyle\qquad\qquad\qquad\qquad + \left(\mathcal Q_0^{LH}(f_{L,1}^{n+1}, f_{H,0}^{n+1})(v^L) + \varepsilon\mathcal Q_{\varepsilon}^{LH}(f_{L,1}^n, f_{H,1}^n)(v^L)\right)\bigg] \phi_1^L\, dv^L \notag \\[2pt]
&\displaystyle\label{M_L1_b}\qquad\qquad\qquad\qquad +\frac{\mu\Delta t}{\varepsilon^2}\Big((P_1)_{L,1}^n - (P_1)_{L,1}^{n+1}\Big),  \\[2pt]
&\displaystyle (P_2)_{L,1}^{n+1} = (P_2)_{L,1}^n + \frac{\Delta t}{\varepsilon^2}\int_{\mathbb R^3}\bigg[\frac{1}{\varepsilon}\, 
\mathcal Q_{\varepsilon}^{LH}(f_{L,0}^{n+1}, f_{H,0}^{n+1})(v^L)
+ \mathcal Q_{\varepsilon}^{LH}(f_{L,0}^n, f_{H,1}^n)(v^L) \notag\\[2pt]
&\displaystyle \qquad\qquad\qquad\qquad\qquad\qquad 
+ \mathcal Q_{\varepsilon}^{LH}(f_{L,1}^n, f_{H,0}^n)(v^L) + \varepsilon\mathcal Q_{\varepsilon}^{LH}(f_{L,1}^n, f_{H,1}^n)(v^L)\bigg] \phi_2^L\, dv^L \notag\\[2pt]
&\displaystyle\label{M_L1_c}\qquad\qquad\qquad\qquad\qquad\qquad +\frac{\mu\Delta t}{\varepsilon^2}\Big((P_2)_{L,1}^n - (P_2)_{L,1}^{n+1}\Big), 
\end{align}
The reason why the BGK- or Fokker-Planck type penalties do not work well for $f_{L,1}$
is due to the complexity of the moment equation (\ref{M_L1_b}), in which
the term $\mathcal Q_0^{LH}(f_{L,1}^{n+1}, f_{H,0}^{n+1})(v^L)$ is implicit since $f_{L,1}^{n+1}$ is unknown. 
We find it difficult to invert this term, since both the moment equation (\ref{M_L1_b}) and the equation (\ref{num2N}) for $f_{L,1}$ involve the same term $f_{L,1}^{n+1}$, thus the entire coupled system (\ref{num1N})--(\ref{num2N}) need to be inverted all together. 
Thus it is hard to get the Maxwellian associated with $f_{L,0}+f_{L,1}$ in the BGK- or Fokker-Planck type penalty operators. 
Investigating a better approach than the currently used linear penalty method in (\ref{LP}) is deferred to a future work. 

For the second collision term $\mathcal Q^{LL}(f_{L,0}^{n}-f_{L,1}^{n}, f_{L,0}^{n}-f_{L,1}^{n})$ in (\ref{Q-LL}), 
the reason we adopt the linear penalty is to avoid negative values of the temperature difference computed from the moment equations of 
$f_{L,0}$ and $f_{L,1}$ (hence unable to define the Maxwellian in the penalty operators). 
The difference between the Filbet-Jin (or Jin-Yan) penalty and the linear penalty is that the latter owns an error of $O(\Delta t)$ compared to $O(\varepsilon)$ as in the former, in the AP analysis. See \cite{Filbet-Jin}. 
Another disadvantage of the linear penalty method is that 
the linear operator does not preserve exactly the mass, momentum and energy as the BGK-type operator does, 
as mentioned in \cite{Filbet-Jin}. Nevertheless, the conservation issues (conservation of mass for each species, and conservation of total momentum and energy for the two species) will be addressed in our follow-up work. 
\end{remark}


\subsubsection{The final numerical scheme}

\indent To summarize, the schemes for $f_{L,0}$, $f_{L,1}$ are given by 
\begin{align}
&\displaystyle\label{num1N-P}\frac{f_{L,0}^{n+1}-f_{L,0}^n}{\Delta t} =
\frac{1}{\varepsilon^2}\bigg[\mathcal Q^{LL}(f_{L,0}^{n}, f_{L,0}^{n}) - \mathcal P(f_{L,0}^n) + \mathcal P(f_{L,0}^{n+1})
+\mathcal Q_0^{LH}(f_{L,0}^{n}, f_{H,0}^{n})\bigg], \\[4pt]
&\displaystyle\frac{f_{L,1}^{n+1}-f_{L,1}^n}{\Delta t} = \frac{1}{\varepsilon^2}\bigg[\frac{1}{\varepsilon}\left(\mathcal Q^{LH}_{\varepsilon}(f_{L,0}^{n+1}, f_{H,0}^{n+1}) -\mathcal Q_0^{LH}(f_{L,0}^{n+1}, f_{H,0}^{n+1})\right) \notag\\[2pt]
&\displaystyle\qquad\qquad\qquad\quad + \frac{1}{2}\bigg[\mathcal Q^{LL}(f_{L,0}^{n}+f_{L,1}^{n}, f_{L,0}^{n}+f_{L,1}^{n}) + \mu (f_{L,0}^n+f_{L,1}^n) - \mu (f_{L,0}^{n+1}+f_{L,1}^{n+1}) \notag\\[2pt]
&\displaystyle \qquad\qquad\qquad\qquad - \left(\mathcal Q^{LL}(f_{L,0}^{n}-f_{L,1}^{n}, f_{L,0}^{n}-f_{L,1}^{n}) + \mu (f_{L,0}^n-f_{L,1}^n)
- \mu (f_{L,0}^{n+1}-f_{L,1}^{n+1})\right)\bigg] \notag\\[2pt]
&\displaystyle\qquad\qquad\qquad\quad  + \varepsilon\mathcal Q^{LL}(f_{L,1}^{n}, f_{L,1}^{n})  + \mathcal Q^{LH}_{\varepsilon}(f_{L,0}^{n}, f_{H,1}^{n}) + \left(\mathcal Q^{LH}_{\varepsilon}(f_{L,1}^{n}, f_{H,0}^{n}) - \mathcal Q_0^{LH}(f_{L,1}^{n}, f_{H,0}^{n})\right) \notag\\[2pt]
&\displaystyle\label{num2N-P}\qquad\qquad\qquad\quad + \mathcal Q_0^{LH}(f_{L,1}^{n+1}, f_{H,0}^{n+1}) +\varepsilon\mathcal Q^{LH}_{\varepsilon}(f_{L,1}^{n}, f_{H,1}^n)\bigg]. 
\end{align}
The schemes for $f_{H,0}$, $f_{H,1}$ are given by
\begin{align}
&\displaystyle\label{num3N-P}\frac{f_{H,0}^{n+1}-f_{H,0}^n}{\Delta t} = 
\frac{1}{\varepsilon}\bigg[\mathcal Q^{HH}(f_{H,0}^{n+1}, f_{H,0}^{n+1}) - \mathcal P(f_{H,0}^n) + \mathcal P(f_{H,0}^{n+1})
+ \mathcal Q_0^{HL}(f_{H,0}^{n}, f_{L,0}^{n})\bigg], \\[4pt]
&\displaystyle\frac{f_{H,1}^{n+1}-f_{H,1}^n}{\Delta t} = \frac{1}{\varepsilon}\bigg[\frac{1}{\varepsilon}\left(\mathcal Q^{HL}_{\varepsilon}(f_{H,0}^{n+1}, f_{L,0}^{n+1})-\mathcal Q_0^{HL}(f_{H,0}^{n+1}, f_{L,0}^{n+1})\right)\notag\\[2pt]
&\displaystyle\qquad\qquad\qquad\quad
+ \frac{1}{2}\bigg[\mathcal Q^{HH}(f_{H,0}^{n}+f_{H,1}^{n}, f_{H,0}^{n}+f_{H,1}^{n}) + \mu(f_{H,0}^n+f_{H,1}^n)  
- \mu(f_{H,0}^{n+1}+f_{H,1}^{n+1}) \notag\\[2pt]
&\displaystyle\qquad\qquad\qquad\qquad  - \left(\mathcal Q^{HH}(f_{H,0}^{n}-f_{H,1}^{n}, f_{H,0}^{n}-f_{H,1}^{n}) + \mu(f_{H,0}^n-f_{H,1}^n)
- \mu(f_{H,0}^{n+1}-f_{H,1}^{n+1})\right)\bigg] \notag \\[2pt]
&\displaystyle\label{num4N-P} \qquad\qquad\qquad\quad  + \varepsilon\mathcal Q^{HH}(f_{H,1}^{n}, f_{H,1}^{n}) + \mathcal Q^{HL}_{\varepsilon}(f_{H,0}^{n+1}, f_{L,1}^{n+1}) + \mathcal Q^{HL}_{\varepsilon}(f_{H,1}^{n}, f_{L,0}^{n}) + \varepsilon\mathcal Q^{HL}_{\varepsilon}(f_{H,1}^{n}, f_{L,1}^n)\bigg]. 
\end{align}
We couple with the following equations for moments of $f_{L,0}$ and $f_{H,0}$ (recall (\ref{Moments}) for the definition): 
\begin{align}
&\displaystyle\label{M_L0_a} (P_0)_{L,0}^{n+1} = (P_0)_{L,0}^n, \\[2pt]
&\displaystyle\label{M_L0_b} (P_1)_{L,0}^{n+1} = (P_1)_{L,0}^n + \frac{\Delta t}{\varepsilon^2}\int_{\mathbb R^3} v^L\, 
 \mathcal Q_0^{LH}(f_{L,0}^n, f_{H,0}^n)(v^L)\, dv^L,  \\[2pt]
&\displaystyle\label{M_L0_c} (P_2)_{L,0}^{n+1} = (P_2)_{L,0}^n, \\[4pt]
&\displaystyle\label{M_H0_a} (P_0)_{H,0}^{n+1} = (P_0)_{H,0}^n, \\[2pt]
&\displaystyle\label{M_H0_b} (P_1)_{H,0}^{n+1} = (P_1)_{H,0}^n + \frac{\Delta t}{\varepsilon^2}\int_{\mathbb R^3}
v^H\, \mathcal Q_0^{HL}(f_{H,0}^n, f_{L,0}^n)(v^H)\, dv^H, \\[2pt]
&\displaystyle\label{M_H0_c} (P_2)_{H,0}^{n+1} = (P_2)_{H,0}^n + \frac{\Delta t}{\varepsilon^2}\int_{\mathbb R^3}
\frac{|v^H|^2}{2}\, \mathcal Q_0^{HL}(f_{H,0}^n, f_{L,0}^n)(v^H)\, dv^H. 
\end{align}
From the moment system, one computes $u$ from $u=\frac{P_1}{P_0}$ and solves for $T$ by using the formula
$$ P_2 = \frac{1}{2} P_0 |u|^2 + \frac{3}{2} P_0\, T, $$
then obtain the local Maxwellian by the definition 
$$ M_{n,u,T}(v)=\frac{n}{(2\pi T)^{3/2}}\,\exp\left(-\frac{|v-u|^2}{2T}\right). $$
$M_{L,0}^{n+1}$ (or $M_{H,0}^{n+1}$) is obtained by $n$, $u$, $T$ 
got from the moment equations of $f_{L,0}$ (or $f_{H,0}$), namely (\ref{M_L0_a})--(\ref{M_L0_c}) (or (\ref{M_H0_a})--(\ref{M_H0_c})). 
\\[2pt]

The following shows the detailed steps for the implementation of our proposed numerical scheme:  \\[2pt]
(a) get $M_{L,0}^{n+1}$ from (\ref{M_L0_a})--(\ref{M_L0_c}), then update $f_{L,0}^{n+1}$ from (\ref{num1N-P});  \\
(b) get $M_{H,0}^{n+1}$ from (\ref{M_H0_a})--(\ref{M_H0_b}), then update $f_{H,0}^{n+1}$ from (\ref{num3N-P});  \\
(c) update $f_{L,1}^{n+1}$ from (\ref{num2N-P}); \\
(d) update $f_{H,1}^{n+1}$ from (\ref{num4N-P}). 
\\[2pt]

Our scheme, although contains some implicit terms, can be implemented {\it explicitly} for the case of Boltzmann collision operator, or just needs a linear elliptic solver in the case of FPL operator, as in the case of single species counterpart
in \cite{Filbet-Jin} and \cite{JinYan}. We would like to mention that higher order time approximation can be extended. 

\subsection{The AP Property}

Our goal of this subsection is to prove the AP property of the discretized scheme (\ref{num1N})--(\ref{num2N}) and (\ref{num3N})--(\ref{num4N}). 

First, for the light particles, inserting the expansion
$$\mathcal Q_{\varepsilon}^{LH}=\mathcal Q_0^{LH} +\varepsilon \mathcal Q_1^{LH} +
\varepsilon^2\, \mathcal Q_2^{LH} +O(\varepsilon^3) $$ 
into (\ref{num2N}), one has
\begin{align}
&\displaystyle\quad\frac{f_{L,1}^{n+1}-f_{L,1}^n}{\Delta t}  \notag\\[4pt]
&\displaystyle = \frac{1}{\varepsilon^2}\Big[\, 2\mathcal Q^{LL}(f_{L,0}^{n+1}, f_{L,1}^{n+1})
+\mathcal Q_0^{LH}(f_{L,0}^n, f_{H,1}^n) + \mathcal Q_0^{LH}(f_{L,1}^{n+1}, f_{H,0}^{n+1}) + \mathcal Q_1^{LH}(f_{L,0}^{n+1}, f_{H,0}^{n+1})\Big]\notag\\[4pt]
&\displaystyle\quad+ \frac{1}{\varepsilon}\Big[\mathcal Q^{LL}(f_{L,1}^n, f_{L,1}^n) +\mathcal Q_0^{LH}(f_{L,1}^n, f_{H,1}^n)
+ \mathcal Q_1^{LH}(f_{L,0}^n, f_{H,1}^n) + \mathcal Q_1^{LH}(f_{L,1}^n, f_{H,0}^n) + \mathcal Q_2^{LH}(f_{L,0}^{n+1}, f_{H,0}^{n+1})\Big] \notag\\[4pt]
&\displaystyle\label{AP_L0}\quad + \mathcal Q_1^{LH}(f_{L,1}^n, f_{H,1}^n) + \mathcal Q_2^{LH}(f_{L,0}^n, f_{H,1}^n)
+\mathcal Q_2^{LH}(f_{L,1}^n, f_{H,0}^n). 
\end{align}
\\[2pt]
First, (\ref{num1N}) gives 
$$\mathcal Q^{LL}(f_{L,0}^{n+1}, f_{L,0}^{n+1})+
\mathcal Q_0^{LH}(f_{L,0}^{n}, f_{H,0}^{n}) = O(\varepsilon^2), $$
thus 
\begin{equation}\label{f_L0} f_{L,0}^{n+1} = n_{L,0}^{n+1}\, M_{0, T_{L,0}^{n+1}} + O(\varepsilon^2 + \Delta t) := M_{L,0}^{n+1} +
O(\varepsilon^2 + \Delta t). \end{equation}
As for the heavy particles, by (\ref{num3N}), 
$$\mathcal Q^{HH}(f_{H,0}^{n+1}, f_{H,0}^{n+1})+
\underbrace{\mathcal Q_0^{HL}(f_{H,0}^{n}, f_{L,0}^{n})}_{=O(\varepsilon^2+\Delta t)}= 
O(\varepsilon), $$
which gives
\begin{equation}\label{f_H0}
f_{H,0}^{n+1} = n_{H,0}^{n+1}\, M_{u_{H,0}^{n+1}, T_{H,0}^{n+1}} + O(\varepsilon+\Delta t) := M_{H,0}^{n+1} +
O(\varepsilon + \Delta t). 
\end{equation}
According to (\ref{AP_L0}), 
\begin{align}
&\displaystyle\label{eps-1} \quad 2\mathcal Q^{LL}(f_{L,0}^{n+1}, f_{L,1}^{n+1})
+ \underbrace{\mathcal Q_0^{LH}(f_{L,0}^{n}, f_{H,1}^{n})}_{=O(\varepsilon^2 + \Delta t)} + \mathcal Q_0^{LH}(f_{L,1}^{n+1}, f_{H,0}^{n+1}) + 
\mathcal Q_1^{LH}(f_{L,0}^{n+1}, f_{H,0}^{n+1}) = O(\varepsilon^2), 
\end{align}
which is an equation for $f_{L,1}^{n+1}$, and can be equivalently written in the form 
$$\phi_{L}^{n+1}=f_{L,1}^{n+1}\, (M_{L,0}^{n+1})^{-1}$$ with 
$$\Gamma_{L,0}\, \phi_{L}^{n+1} = -(M_{L,0}^{n+1})^{-1}\, \mathcal Q_1^{LH}(M_{L,0}^{n+1}, \, M_{H,0}^{n+1}) + 
O(\varepsilon + \Delta t), $$
where $\Gamma_{L,0}$ is the linearized operator
$$\Gamma_{L,0}\, \phi_{L}^{n+1} = (M_{L,0}^{n+1})^{-1} \left[ 2\mathcal Q^{LL}(M_{L,0}^{n+1},\, M_{L,0}^{n+1}\, \phi_{L}^{n+1}) + 
\mathcal Q_0^{LH}(M_{L,0}^{n+1}\, \phi_{L}^{n+1},\, M_{H,0}^{n+1})\right]. $$ 
Analogous to the continuous case proved in \cite{Degond}, the unique solution in $\left(\text{ker}(\Gamma_{L,0})\right)^{\perp}$ is given by
\begin{equation}
\label{f_L1} f_{L,1}^{n+1}(v^L) = \underbrace{\frac{M_{L,0}^{n+1}}{T_{L,0}^{n+1}(v^L)}\, u_{H,0}^{n+1}\cdot v^L}_{:=f_{L,1}^{\ast, n+1}}\, +\, 
O(\varepsilon + \Delta t), \end{equation}
where $f_{L,1}^{\ast, n+1}$ is used to denote the leading order of $f_{L,1}^{n+1}$. 
\\[2pt]

Multiply (\ref{AP_L0}) by $\varepsilon$ and add up with (\ref{num1N}), then 
\begin{align}
&\displaystyle\quad \frac{f_{L,0}^{n+1}-f_{L,0}^n}{\Delta t} +\varepsilon\, \frac{f_{L,1}^{n+1}-f_{L,1}^n}{\Delta t}  \notag\\[4pt]
&\displaystyle = \frac{1}{\varepsilon^2}\Big[\mathcal Q^{LL}(f_{L,0}^{n+1}, f_{L,0}^{n+1})+
\mathcal Q_0^{LH}(f_{L,0}^{n}, f_{H,0}^{n})\Big] \notag\\[4pt]
&\displaystyle\quad + \frac{1}{\varepsilon}\Big[ 2\mathcal Q^{LL}(f_{L,0}^{n+1}, f_{L,1}^{n+1})
+\mathcal Q_0^{LH}(f_{L,0}^n, f_{H,1}^n) + \mathcal Q_0^{LH}(f_{L,1}^{n+1}, f_{H,0}^{n+1}) + \mathcal Q_1^{LH}(f_{L,0}^{n+1}, f_{H,0}^{n+1})\Big]\notag\\[4pt]
&\displaystyle\quad + \mathcal Q^{LL}(f_{L,1}^n, f_{L,1}^n) +\mathcal Q_0^{LH}(f_{L,1}^n, f_{H,1}^n)
+ \mathcal Q_1^{LH}(f_{L,0}^n, f_{H,1}^n) + \mathcal Q_1^{LH}(f_{L,1}^n, f_{H,0}^n) + \mathcal Q_2^{LH}(f_{L,0}^{n+1}, f_{H,0}^{n+1}) \notag\\[4pt]
&\displaystyle\label{AP_L}\quad + \varepsilon\left[\mathcal Q_1^{LH}(f_{L,1}^n, f_{H,1}^n) + \mathcal Q_2^{LH}(f_{L,0}^n, f_{H,1}^n)
+\mathcal Q_2^{LH}(f_{L,1}^n, f_{H,0}^n)\right]. 
\end{align}
Plugging in the leading order of (\ref{f_L0}), (\ref{f_L1}) and comparing the $O(1)$ terms on both sides, one gets 
\begin{align}
&\displaystyle \frac{M_{L,0}^{n+1}-M_{L,0}^n}{\Delta t} 
= \mathcal Q^{LL}(f_{L,1}^{\ast, n}, f_{L,1}^{\ast, n}) + \mathcal Q_0^{LH}(f_{L,1}^{\ast, n}, f_{H,1}^{\ast, n})
+ \mathcal Q_1^{LH}(M_{L,0}^n, f_{H,1}^{\ast, n}) \notag\\[4pt]
&\displaystyle\label{eps0}\qquad\qquad\qquad\qquad + \mathcal Q_1^{LH}(f_{L,1}^{\ast, n}, M_{H,0}^n) +  \mathcal Q_2^{LH}(M_{L,0}^{n+1}, M_{H,0}^{n+1}) 
+ O(\Delta t). 
\end{align}
Integrate both sides of (\ref{eps0}) against $|v^L|^2$ on $v^L$, then 
\begin{align*}
&\displaystyle \int \mathcal Q^{LL}(f_{L,1}^{\ast, n}, f_{L,1}^{\ast, n})\, |v^L|^2\, dv^L = \int \mathcal Q_0^{LH}(f_{L,1}^{\ast, n}, f_{H,1}^{\ast, n})\, |v^L|^2\, dv^L = 0, \\[4pt]
&\displaystyle  \int  \mathcal Q_1^{LH}(M_{L,0}^n, f_{H,1}^{\ast, n})\, |v^L|^2\, dv^L = \int \mathcal Q_0^{HL}(f_{H,1}^{\ast, n}, M_{L,0}^n)\, |v^L|^2\, dv^L =0, \end{align*}
and 
\begin{align*}
&\displaystyle \quad\int\left[\mathcal Q_1^{LH}(f_{L,1}^{\ast, n}, M_{H,0}^{n}) +  \mathcal Q_2^{LH}(M_{L,0}^{n+1}, M_{H,0}^{n+1})\right] |v^L|^2\, dv^L  \\[4pt]
&\displaystyle = \int \left[\mathcal Q_0^{HL}(M_{H,0}^{n}, f_{L,1}^{\ast, n}) + \mathcal Q_1^{HL}(M_{H,0}^{n+1}, M_{L,0}^{n+1})\right] |v^H|^2\, dv^H \\[4pt]
&\displaystyle = \int \left[\mathcal Q_0^{HL}(M_{H,0}^{n+1}, f_{L,1}^{\ast, n+1}) + \mathcal Q_1^{HL}(M_{H,0}^{n+1}, M_{L,0}^{n+1})\right] |v^H|^2\, dv^H 
+ O(\Delta t) \\[4pt]
&\displaystyle = 3\, \frac{\lambda(T_{L,0}^{n+1})}{T_{L,0}^{n+1}}\, n_{L,0}^{n+1}\, n_{H,0}^{n+1}\, (T_{H,0}^{n+1}- T_{L,0}^{n+1}) + O(\Delta t), 
\end{align*}
where analogous calculation of the integrals for the continuous case is shown in \cite{Degond}. 
Denote $\mathcal D_t(u^n)$ the discrete time derivative of the numerical quantity of interest $u^n$: 
$$ \mathcal D_t(u^n) := \frac{u^{n+1} - u^n}{\Delta t}\,. $$

Integrating both sides of (\ref{eps0}) on $v^L$ gives
$$\mathcal D_t (n_{L,0}^n) = O(\Delta t), $$
by using (\ref{Q_i}) in Theorem \ref{property}. 
Integrals of the left-hand-side of (\ref{eps0}) against $1$ and $|v^L|^2$ on $v^L$ are
$$\mathcal D_t \left(n_{L,0}^{n}, \, n_{L,0}^{n}(\frac{1}{2}|u_{L,0}^{n}|^2 +\frac{3}{2}T_{L,0}^{n})\right)^{T}. $$

Therefore, the limit of our discretized numerical scheme is given by 
\begin{align*}
&\displaystyle \mathcal D_t (n_{L,0}^n) = O(\Delta t),  \\[2pt]
&\displaystyle\mathcal D_t\left(\frac{3}{2}n_{L,0}^n\, T_{L,0}^n\right) = 3\, \frac{\lambda(T_{L,0}^{n+1})}{T_{L,0}^{n+1}}\, n_{L,0}^{n+1}\, n_{H,0}^{n+1}\, (T_{H,0}^{n+1}- T_{L,0}^{n+1}) + O(\Delta t), 
\end{align*}
which is consistent with the implicit discretization of the continuous limit (\ref{T2}), up to a numerical error of $O(\Delta t)$. 
\\[2pt]

Now we examine the system for the heavy particles $f_{H,0}$, $f_{H,1}$. 
Multiplying (\ref{num4N}) by $\varepsilon$, adding it up with (\ref{num3N}) and using the expansion
$$\mathcal Q_{\varepsilon}^{HL} = \mathcal Q_0^{HL} + \varepsilon\mathcal Q_1^{HL} + \varepsilon^2\, \mathcal Q_2^{HL} 
+ O(\varepsilon^3), $$ 
one gets
\begin{align}
&\displaystyle\quad \frac{f_{H,0}^{n+1}-f_{H,0}^n}{\Delta t} + \varepsilon\, \frac{f_{H,1}^{n+1}-f_{H,1}^n}{\Delta t}  \notag\\[4pt]
&\displaystyle =  \frac{1}{\varepsilon}\bigg[\mathcal Q^{HH}(f_{H,0}^{n+1}, f_{H,0}^{n+1})+
\mathcal Q_0^{HL}(f_{H,0}^{n}, f_{L,0}^{n})\bigg] \notag\\[4pt]
&\displaystyle\quad + \mathcal Q_0^{HL}(f_{H,0}^{n+1}, f_{L,1}^{n+1}) + \mathcal Q_0^{HL}(f_{H,1}^n, f_{L,0}^n) 
+ 2\mathcal Q^{HH}(f_{H,0}^{n+1}, f_{H,1}^{n+1}) + \mathcal Q_1^{HL}(f_{H,0}^{n+1}, f_{L,0}^{n+1}) \notag\\[4pt]
&\displaystyle\quad + \varepsilon\bigg[\mathcal Q^{HH}(f_{H,1}^n, f_{H,1}^n) +\mathcal Q_1^{HL}(f_{H,0}^{n+1}, f_{L,1}^{n+1})
+\mathcal Q_1^{HL}(f_{H,1}^n, f_{L,0}^n) + \mathcal Q_0^{HL}(f_{H,1}^n, f_{L,1}^n) + \mathcal Q_2^{HL}(f_{H,0}^{n+1}, f_{L,0}^{n+1})\bigg] \notag\\[4pt]
&\displaystyle\label{AP_H}\quad + \varepsilon^2\left[\mathcal Q_1^{HL}(f_{H,1}^n, f_{L,1}^n) + \mathcal Q_2^{HL}(f_{H,0}^{n+1}, f_{L,1}^{n+1}) + 
\mathcal Q_2^{HL}(f_{H,1}^{n+1}, f_{L,0}^{n+1})\right]. 
\end{align}
\\[4pt]
Plug in the leading order term of (\ref{f_H0}) and compare the $O(1)$ terms on both sides, then 
\begin{align*}
&\displaystyle \frac{M_{H,0}^{n+1}-M_{H,0}^n}{\Delta t} = 
 2\mathcal Q^{HH}(M_{H,0}^{n+1}, f_{H,1}^{\ast, n+1}) + \mathcal Q_0^{HL}(M_{H,0}^{n+1}, f_{L,1}^{\ast, n+1}) + 
 \underbrace{\mathcal Q_0^{HL}(f_{H,1}^{\ast, n}, M_{L,0}^{n})}_{=0}
 \\[4pt]
&\displaystyle\qquad\qquad\qquad\quad +\mathcal Q_1^{HL}(M_{H,0}^{n+1}, M_{L,0}^{n+1}) +  O(\Delta t). 
\end{align*}
It is an equation for $f_{H,1}^{\ast, n+1}$ and can be equivalently written in terms of 
$$\phi_{H}^{n+1} = f_{H,1}^{\ast, n+1}\, (M_{H,0}^{n+1})^{-1}$$
according to 
\begin{equation}\label{Gamma_H0}\Gamma_{H,0}\, \phi_{H}^{n+1} = (M_{H,0}^{n+1})^{-1}\left[\mathcal D_t M_{H,0}^{n} - \mathcal Q_0^{HL}(M_{H,0}^{n+1}, \, f_{L,1}^{\ast, n+1}) - \mathcal Q_1^{HL}(M_{H,0}^{n+1}, M_{L,0}^{n+1})\right] +  O(\Delta t), 
 \end{equation}
$\Gamma_{H,0}$ is a linearization operator given by 
$$\Gamma_{H,0}\, \phi_{H}^{n+1} = 2\, (M_{H,0}^{n+1})^{-1}\, \mathcal Q^{HH}(M_{H,0}^{n+1}, \, M_{H,0}^{n+1}\, \phi_{H}^{n+1})\,.$$
The necessary and sufficient condition for the solvability of equation (\ref{Gamma_H0}) is given by
\begin{equation}
\label{H0_cond1}\int_{\mathbb R^3}\,\left[D_t M_{H,0}^{n} - \mathcal Q_0^{HL}(M_{H,0}^{n+1}, \, f_{L,1}^{\ast, n+1}) - \mathcal Q_1^{HL}(M_{H,0}^{n+1}, M_{L,0}^{n+1})\right]\begin{pmatrix} 1\\ v^H \\|v^H|^2\end{pmatrix} dv^H = O(\Delta t)\, \mathbb I_{3}, 
 \end{equation}
 where $\mathbb I_{3}=\left(1, 1, 1\right)^{T}$. 
Analogous to the calculation in \cite{Degond} for the continuous equations, 
\begin{align}
&\displaystyle \,\int_{\mathbb R^3}\, \left[\mathcal Q_0^{HL}(M_{H,0}^{n+1}, \, f_{L,1}^{\ast, n+1}) +\mathcal Q_1^{HL}(M_{H,0}^{n+1}, M_{L,0}^{n+1})\right](v^H)
\begin{pmatrix}1\\v^H\\ \frac{|v^H|^2}{2}\end{pmatrix}dv^H \notag\\[4pt]
&\displaystyle\label{H0_cond2}  = \begin{pmatrix} 0 \\ 0 \\ -3\,\frac{\lambda(T_{L,0}^{n+1})}{T_{L,0}^{n+1}}\, n_{L,0}^{n+1} n_{H,0}^{n+1}(T_{H,0}^{n+1}- T_{L,0}^{n+1})\end{pmatrix} +  O(\Delta t)\, \mathbb I_{3}. 
\end{align}
Insert (\ref{H0_cond2}) into (\ref{H0_cond1}), then
\begin{equation}
\label{Limit_H}
\mathcal D_t \begin{pmatrix} n_{H,0}^{n}\\ n_{H,0}^{n}\, u_{H,0}^{n}\\ n_{H,0}^{n}(\frac{1}{2}|u_{H,0}^{n}|^2 +\frac{3}{2}T_{H,0}^{n})
\end{pmatrix}= \begin{pmatrix} 0\\0\\ -3\,\frac{\lambda(T_{L,0}^{n+1})}{T_{L,0}^{n+1}}\, n_{L,0}^{n+1} n_{H,0}^{n+1}(T_{H,0}^{n+1}- T_{L,0}^{n+1})\end{pmatrix}
+  O(\Delta t)\, \mathbb I_{3}. 
\end{equation}
This shows that $n_{H,0}^{n}$, $u_{H,0}^{n}$ are constant in time with a numerical error of $O(\Delta t)$, 
$$\mathcal D_t(n_{H,0}^n) = O(\Delta t),  \qquad \mathcal D_t(n_{H,0}^n u_{H,0}^n) = O(\Delta t), $$
while $T_{H,0}^n$ evolves according to 
$$\mathcal D_t \left(\frac{3}{2}n_{H,0}^{n}T_{H,0}^{n}\right) = -3\,\frac{\lambda(T_{L,0}^{n+1})}{T_{L,0}^{n+1}}\, n_{L,0}^{n+1} n_{H,0}^{n+1}(T_{H,0}^{n+1}- T_{L,0}^{n+1}) + O(\Delta t), $$
which is consistent with the discretized implicit scheme of the limiting system (\ref{T1}), up to a numerical error of $O(\Delta t)$. 

We conclude our AP analysis with the following theorem. 
\begin{theorem}
\label{main}
The time discretized numerical schemes given by 
(\ref{num1N})--(\ref{num2N}) and (\ref{num3N})--(\ref{num4N}), 
as $\varepsilon\to 0$, approaches to the system 
\begin{align*}
&\displaystyle n_{L,0}^{n+1}=n_{L,0}^n + O(\Delta t), \\[2pt]
&\displaystyle n_{H,0}^{n+1}=n_{H,0}^n + O(\Delta t), \qquad n_{H,0}^{n+1}\, u_{H,0}^{n+1}=n_{H,0}^n\, u_{H,0}^n + O(\Delta t), \\[2pt]
&\displaystyle \frac{d}{dt}\left(\frac{3}{2}n_{L,0}^n\, T_{L,0}^n\right) = 3\, \frac{\lambda(T_{L,0}^{n+1})}{T_{L,0}^{n+1}}\, n_{L,0}^{n+1}\, n_{H,0}^{n+1}\, (T_{H,0}^{n+1}- T_{L,0}^{n+1}) + O(\Delta t), \\[2pt]
&\displaystyle\frac{d}{dt}\left(\frac{3}{2}n_{H,0}^{n}T_{H,0}^{n}\right) = -3\,\frac{\lambda(T_{L,0}^{n+1})}{T_{L,0}^{n+1}}\, n_{L,0}^{n+1} n_{H,0}^{n+1}
(T_{H,0}^{n+1}- T_{L,0}^{n+1}) + O(\Delta t), 
\end{align*}
which are consistent with the implicit discretization of the continuous limit shown in (\ref{T2})--(\ref{T1}), 
with a numerical error of $O(\Delta t)$. 
\end{theorem}

\begin{remark}
We would also like to point out that our AP analysis for the scheme does not include the penalty method, namely the schemes (\ref{num1N-P})--(\ref{num4N-P}) that one actually uses in practice, 
since it is hard to prove a scheme is AP with all the penalty terms included, not done even for the single species Boltzmann (or FPL)
equation \cite{Filbet-Jin, JinYan}.  
\end{remark}

\section{The Space Inhomogeneous Systems}
\label{sec3}
In the space inhomogeneous case, the evolution equations are given by system (\ref{SI_L})--(\ref{SI_H}). 
We first recall the main results in \cite{Degond2} in the following Theorem: 
\begin{theorem}
As $\varepsilon\to 0$, the limit distributions and limit systems are given by 
$$ f_0^L(x,v,t) = n_0^L(x,t)\, M_{0, T_0^L(x,t)}, \qquad
f_0^H(x,v,t) = n_0^H(x,t)\, M_{u_0^H(x,t), T_0^H(x,t)}, $$
where $n_0^L$, $T_0^L$, $n_0^H$, $T_0^H$ satisfy the coupled system: 
\begin{align}
&\displaystyle\label{nL0} \frac{\partial n_0^L}{\partial t} + \nabla_x \cdot (n_0^L u_0^H)
- \nabla_x \cdot \left[D_{11}\left(\nabla_x n_0^L - \frac{F^L n_0^L}{T_0^L}\right) + D_{12}\left(n_0^L\, \frac{\nabla_x T_0^L}{T_0^L}\right)\right] = 0, \\[4pt]
&\displaystyle\notag \frac{\partial}{\partial t}\left(\frac{3}{2}n_0^L T_0^L\right)  + \nabla_x \cdot \left(\frac{5}{2}n_0^L T_0^L u_0^H\right) - n_0^L F^L\cdot u_0^H \\[2pt]
&\displaystyle\notag\qquad\qquad\qquad - \nabla_x\cdot\left[D_{21}\left(\nabla_x n_0^L - \frac{F^L n_0^L}{T_0^L}\right) + D_{22}\left(n_0^L\, \frac{\nabla_x T_0^L}{T_0^L}  \right)\right] \\[2pt]
&\displaystyle\notag\qquad\qquad\qquad + F^L \cdot\left[D_{11}\left(\nabla_x n_0^L - \frac{F^L n_0^L}{T_0^L}\right) + D_{12}\left(n_0^L\, \frac{\nabla_x T_0^L}{T_0^L}  \right)\right] \\[2pt]
&\displaystyle\label{nL1} = u_0^H \cdot \left[\nabla_x(n_0^L T_0^L) - F^L n_0^L\right] + 3 \frac{\lambda(T_0^L)}{T_0^L}\, n_0^L n_0^H\, (T_0^H - T_0^L), 
\end{align}
and
\begin{align}
&\displaystyle\label{nH0} \frac{\partial n_0^H}{\partial t} + \nabla_x \cdot (n_0^H u_0^H) = 0, \\[4pt]
&\displaystyle\notag \frac{\partial}{\partial t}(n_0^H u_0^H) + \nabla_x \cdot (n_0^H u_0^H\otimes u_0^H)
+ \nabla_x (n_0^H T_0^H) - n_0^H F^H \\[2pt]
&\displaystyle\label{nH1}= - \left(\nabla_x(n_0^L T_0^L) - F^L n_0^L\right), \\[4pt]
&\displaystyle\notag \frac{\partial}{\partial t} \left(\frac{n_0^H |u_0^H|^2}{2} + \frac{3}{2}n_0^H T_0^H\right) 
+ \nabla_x \cdot\left(\left(\frac{n_0^H |u_0^H|^2}{2} + \frac{5}{2}n_0^H T_0^H \right)u_0^H\right) - n_0^H F^H\cdot u_0^H \\[2pt]
&\displaystyle\label{nH2} = -u_0^H \cdot\left[\nabla_x(n_0^L T_0^L) - F^L n_0^L\right] 
- 3 \frac{\lambda(T_0^L)}{T_0^L}\, n_0^L n_0^H\, (T_0^H - T_0^L), 
\end{align}
where $D_{ij}$ ($i,j=1,2$) and $\lambda(T)$ are given in the Appendix. \\[2pt]
\end{theorem}

Insert the expansion $$f^L=f_0^L +\varepsilon f_1^L\,, \qquad f^H=f_0^H+\varepsilon f_1^H$$ into (\ref{SI_L}) and (\ref{SI_H}), then 
\begin{align*}
&\displaystyle\quad \frac{\partial}{\partial t}(f_0^L+\varepsilon f_1^L) + \frac{1}{\varepsilon}\left(v^L\cdot\nabla_x f_0^L + F^L\cdot\nabla_{v^L}f_0^L\right)
+ \left(v^L\cdot\nabla_x f_1^L + F^L\cdot\nabla_{v^L}f_1^L\right) \\[4pt]
&\displaystyle =\frac{1}{\varepsilon^2}\left[\mathcal Q^{LL}(f_0^L+\varepsilon f_1^L, f_0^L+\varepsilon f_1^L)   
+\mathcal Q^{LH}_{\varepsilon}(f_0^L+\varepsilon f_1^L, f_0^H+\varepsilon f_1^H)\right] \\[4pt]
&\displaystyle  = \frac{1}{\varepsilon^2}\Big[\mathcal Q^{LL}(f_0^L, f_0^L) +2\varepsilon \mathcal Q^{LL}(f_0^L, f_1^L) +\varepsilon^2 \mathcal Q^{LL}(f_1^L, f_1^L)\\[4pt]
&\displaystyle\quad +\mathcal Q^{LH}_{\varepsilon}(f_0^L, f_0^H) +\varepsilon \mathcal Q^{LH}_{\varepsilon}(f_0^L, f_1^H) +\varepsilon \mathcal Q^{LH}_{\varepsilon}(f_1^L, f_0^H)
+\varepsilon^2 \mathcal Q^{LH}_{\varepsilon}(f_1^L, f_1^H)\Big]. 
\end{align*}

We design the scheme by letting $f_0^L$, $f_1^L$ satisfy the system 
\begin{align}
&\displaystyle\label{SI_L0}\quad \frac{\partial}{\partial t}f_0^L + \left(v^L\cdot\nabla_x f_1^L + F^L\cdot\nabla_{v^L}f_1^L \right) =\frac{1}{\varepsilon^2}\Big[\mathcal Q^{LL}(f_0^L, f_0^L)+\mathcal Q_0^{LH}(f_0^L, f_0^H)\Big], \\[6pt]
&\displaystyle\quad \frac{\partial}{\partial t}f_1^L + \frac{1}{\varepsilon^2}\left(v^L\cdot\nabla_x f_0^L + F^L\cdot\nabla_{v^L}f_0^L\right) \notag\\[4pt]
 &\displaystyle  = \frac{1}{\varepsilon^2}\Big[\frac{1}{\varepsilon}\left(\mathcal Q^{LH}_{\varepsilon}(f_0^L, f_0^H)-\mathcal Q_0^{LH}(f_0^L, f_0^H)   \right) +2\mathcal Q^{LL}(f_0^L, f_1^L)+\varepsilon\mathcal Q^{LL}(f_1^L, f_1^L)  \notag\\[4pt]
&\displaystyle \label{SI_L1}\qquad +\mathcal Q^{LH}_{\varepsilon}(f_0^L, f_1^H) +
\mathcal Q^{LH}_{\varepsilon}(f_1^L, f_0^H) +\varepsilon\mathcal Q^{LH}_{\varepsilon}(f_1^L, f_1^H)\Big], 
\end{align}
and letting $f_0^H$, $f_1^H$ satisfy the following system
\begin{align} 
&\displaystyle\label{SI_H0}\quad \frac{\partial}{\partial t}f_0^H + \varepsilon\left(v^H\cdot\nabla_x f_1^H + F^H\cdot\nabla_{v^H}f_1^H \right)
=\frac{1}{\varepsilon}\Big[\mathcal Q^{H}(f_0^H, f_0^H)+\mathcal Q_0^{HL}(f_0^H, f_0^L)\Big], \\[6pt]
&\displaystyle\quad \frac{\partial}{\partial t}f_1^H +\frac{1}{\varepsilon}\left(v^H\cdot\nabla_x f_0^H + F^H\cdot\nabla_{v^H}f_0^H \right) \notag\\[4pt]
&\displaystyle = \frac{1}{\varepsilon}\Big[\frac{1}{\varepsilon}\left(\mathcal Q^{HL}_{\varepsilon}(f_0^H, f_0^L)-\mathcal Q_0^{HL}(f_0^H, f_0^L)\right) 
+ 2\mathcal Q^{HH}(f_0^H, f_1^H)+\varepsilon \mathcal Q^{HH}(f_1^H, f_1^H) \notag\\[4pt]
&\displaystyle\label{SI_H1}\quad +\mathcal Q^{HL}_{\varepsilon}(f_0^H, f_1^L) + \mathcal Q^{HL}_{\varepsilon}(f_1^H, f_0^L) + \varepsilon\mathcal Q^{HL}_{\varepsilon}(f_1^H, f_1^L)\Big].  
\end{align}

\subsection{Time discretization}

Following \cite{JPT2}, we rewrite (\ref{SI_L1}) into the diffusive relaxation system
\begin{align}
&\displaystyle \quad\frac{\partial}{\partial t}f_1^L + \psi_1\left(v^L\cdot\nabla_x f_0^L + F^L\cdot\nabla_{v^L}f_0^L\right) \notag\\[4pt]
 &\displaystyle = \frac{1}{\varepsilon^2}\Big[\frac{1}{\varepsilon}\left(\mathcal Q^{LH}_{\varepsilon}(f_0^L, f_0^H)-\mathcal Q_0^{LH}(f_0^L, f_0^H)\right) + 2\mathcal Q^{LL}(f_0^L, f_1^L)+\varepsilon\mathcal Q^{LL}(f_1^L, f_1^L) \notag\\[4pt]
&\displaystyle\qquad +\mathcal Q^{LH}_{\varepsilon}(f_0^L, f_1^H) +\mathcal Q^{LH}_{\varepsilon}(f_1^L, f_0^H) +\varepsilon\mathcal Q^{LH}_{\varepsilon}(f_1^L, f_1^H) \notag\\[4pt]
&\displaystyle \label{relaxation}\qquad -(1-\varepsilon^2\psi_1)\left(v^L\cdot\nabla_x f_0^L + F^L\cdot\nabla_{v^L}f_0^L\right)\Big], 
\end{align}
where a simple choice of $\psi_1$ is $$\psi_1=\min\{1, \, \frac{1}{\varepsilon^2}\}. $$
Note that when $\varepsilon$ is small, $\psi_1=1$.  
The collision operators on the right-hand-side are discretized exactly the same as the space homogeneous case. Then the time discretization
for (\ref{SI_L0}) and (\ref{relaxation}) are 
\begin{align}
&\displaystyle\quad \frac{f_{L,0}^{n+1}-f_{L,0}^n}{\Delta t} + \big(v^L\cdot\nabla_x f_{L,1}^{n} + F_{L}^{n}\cdot\nabla_{v^L}f_{L,1}^{n}\big)
\notag\\[4pt]
&\displaystyle\label{SI_num1N}=\frac{1}{\varepsilon^2}\big[\mathcal Q^{LL}(f_{L}^{n+1}, f_{L,0}^{n+1})+
\mathcal Q_0^{LH}(f_{L,0}^{n}, f_{H,0}^{n})\big]\,, \\[10pt]
&\displaystyle\quad \frac{f_{L,1}^{n+1}-f_{L,1}^n}{\Delta t} + \psi_1\left(v^L\cdot\nabla_x f_{L,0}^{n} + F_{L}^{n}\cdot\nabla_{v^L}f_{L,0}^{n}\right) \notag\\[4pt]
&\displaystyle = \frac{1}{\varepsilon^2}\Big[\frac{1}{\varepsilon}\big(\mathcal Q^{LH}_{\varepsilon}(f_{L,0}^{n+1}, f_{H,0}^{n+1})
-\mathcal Q_0^{LH}(f_{L,0}^{n+1}, f_{H,0}^{n+1})\Big) + 2\mathcal Q^{LL}(f_{L,0}^{n+1}, f_{L,1}^{n+1})+\varepsilon\mathcal Q^{LL}(f_{L,1}^{n}, f_{L,1}^{n}) \notag\\[4pt]
&\displaystyle\qquad\quad + \mathcal Q^{LH}_{\varepsilon}(f_{L,0}^{n}, f_{H,1}^{n}) +
\mathcal Q^{LH}_{\varepsilon}(f_{L,1}^{n+1}, f_{H,0}^{n+1}) + \varepsilon\mathcal Q^{LH}_{\varepsilon}(f_{L,1}^{n}, f_{H,1}^n) \notag\\[4pt]
&\displaystyle\label{SI_num2N}\qquad\quad -(1-\varepsilon^2\psi_1)(v^L\cdot\nabla_x f_{L,0}^{n+1} + F_{L}^{n+1}\cdot\nabla_{v^L}f_{L,0}^{n+1})\Big]. 
\end{align}
\\[2pt]
Using the same technique, time discretization for the systems (\ref{SI_H0}) and (\ref{SI_H1}) are given by 
\begin{align}
&\displaystyle\quad\frac{f_{H,0}^{n+1}-f_{H,0}^n}{\Delta t} +\varepsilon\Big(v^H\cdot\nabla_x f_{H,1}^{n} + 
F_{H}^{n}\cdot\nabla_{v^H}f_{H,1}^{n}\Big)\notag\\[4pt]
&\displaystyle\label{SI_num3N} = \frac{1}{\varepsilon}\Big[\mathcal Q^{HH}(f_{H,0}^{n+1}, f_{H,0}^{n+1})+
\mathcal Q_0^{HL}(f_{H,0}^{n}, f_{L,0}^{n})\Big], \\[10pt]
&\displaystyle\quad\frac{f_{H,1}^{n+1}-f_{H,1}^n}{\Delta t} + {\color{black}\psi_2} \Big(v^H\cdot\nabla_x f_{H,0}^{n} + 
F_{H}^{n}\cdot\nabla_{v^H}f_{H,0}^{n}\Big) \notag\\[4pt]
&\displaystyle = \frac{1}{\varepsilon}\Big[\frac{1}{\varepsilon}\big(\mathcal Q^{HL}_{\varepsilon}(f_{H,0}^{n+1}, f_{L,0}^{n+1})
-\mathcal Q_0^{HL}(f_{H,0}^{n+1}, f_{L,0}^{n+1})\Big)
+ 2\mathcal Q^{HH}(f_{H,0}^{n+1}, f_{H,1}^{n+1}) +\varepsilon\mathcal Q^{HH}(f_{H,1}^{n}, f_{H,1}^{n})\notag \\[4pt]
&\displaystyle\qquad +\mathcal Q^{HL}_{\varepsilon}(f_{H,0}^{n+1}, f_{L,1}^{n+1})+
\mathcal Q^{HL}_{\varepsilon}(f_{H,1}^{n}, f_{L,0}^{n}) + \varepsilon\mathcal Q^{HL}_{\varepsilon}(f_{H,1}^{n}, f_{L,1}^n) \notag \\[4pt]
&\displaystyle\label{SI_num4N}\qquad  {\color{black} - (1-\varepsilon\psi_2)\Big(v^H\cdot\nabla_x f_{H,0}^{n+1} + 
F_{H}^{n}\cdot\nabla_{v^H}f_{H,0}^{n+1}\Big)\Big] }, 
\end{align}
{\color{black}where $$\psi_2=\min\{1, \frac{1}{\varepsilon}\}. $$}
We will use the penalties exactly the same as discussed in subsection \ref{Time-D}, namely the right-hand-side of the 
schemes (\ref{num1N-P})--(\ref{num4N-P}). We omit repeating it here. 

\subsection{The AP Property}

\indent First, for the light particles, inserting the expansion
$$\mathcal Q_{\varepsilon}^{LH} = \mathcal Q_0^{LH} + \varepsilon\mathcal Q_1^{LH} + \varepsilon^2\, \mathcal Q_2^{LH} +
O(\varepsilon^3) $$ 
into (\ref{SI_num2N}), one has
\begin{align}
&\displaystyle \quad \frac{f_{L,1}^{n+1}-f_{L,1}^n}{\Delta t} + 
\big(v^L\cdot\nabla_x f_{L,0}^{n} + F_{L}^{n}\cdot\nabla_{v^L}f_{L,0}^{n}\big)\notag\\[4pt]
&\displaystyle = \frac{1}{\varepsilon^2} \Big[2\mathcal Q^{LL}(f_{L,0}^{n+1}, f_{L,1}^{n+1})
+\mathcal Q_0^{LH}(f_{L,0}^n, f_{H,1}^n) + \mathcal Q_0^{LH}(f_{L,1}^{n+1}, f_{H,0}^{n+1}) + \mathcal Q_1^{LH}(f_{L,0}^{n+1}, f_{H,0}^{n+1}) \notag\\[4pt]
&\displaystyle\qquad - (v^L\cdot\nabla_x f_{L,0}^{n+1} + F_{L}^{n+1}\cdot\nabla_{v^L}f_{L,0}^{n+1})\Big] \notag\\[4pt]
&\displaystyle\qquad + \frac{1}{\varepsilon} \Big[\mathcal Q^{LL}(f_{L,1}^n, f_{L,1}^n) +\mathcal Q_0^{LH}(f_{L,1}^n, f_{H,1}^n)
+ \mathcal Q_1^{LH}(f_{L,0}^n, f_{H,1}^n) + \mathcal Q_1^{LH}(f_{L,1}^n, f_{H,0}^n) + \mathcal Q_2^{LH}(f_{L,0}^{n+1}, f_{H,0}^{n+1})\Big] \notag\\[4pt]
&\displaystyle\label{SI_APL}\qquad + \mathcal Q_1^{LH}(f_{L,1}^n, f_{H,1}^n) + \mathcal Q_2^{LH}(f_{L,0}^n, f_{H,1}^n) 
+ \mathcal Q_2^{LH}(f_{L,1}^n, f_{H,0}^n) + \big(v^L\cdot\nabla_x f_{L,0}^{n+1} + F_{L}^{n+1}\cdot\nabla_{v^L}f_{L,0}^{n+1}\big). 
\end{align}
From (\ref{SI_num1N}), we have 
$$\mathcal Q^{LL}(f_{L,0}^{n+1}, f_{L,0}^{n+1})+
\mathcal Q_0^{LH}(f_{L,0}^{n}, f_{H,0}^{n}) = O(\varepsilon^2), $$
which gives  
\begin{equation}\label{SI_fL0} f_{L,0}^{n+1} = n_{L,0}^{n+1}\, M_{0, T_{L,0}^{n+1}} + O(\varepsilon^2 + \Delta t): = M_{L,0}^{n+1} + 
O(\varepsilon^2 + \Delta t). 
\end{equation}
\\[2pt]
From (\ref{SI_num3N}), $$\mathcal Q^{HH}(f_{H,0}^{n+1}, f_{H,0}^{n+1})+
\underbrace{\mathcal Q_0^{HL}(f_{H,0}^{n}, f_{L,0}^{n})}_{= O(\varepsilon^2 + \Delta t)}= O(\varepsilon), $$
thus \begin{equation}\label{SI_fH0}f_{H,0}^{n+1} = n_{H,0}^{n+1}\, M_{u_{H,0}^{n+1}, T_{H,0}^{n+1}} + O(\Delta t) := M_{H,0}^{n+1} + O(\varepsilon + \Delta t). 
\end{equation}
From (\ref{SI_APL}), 
\begin{align}
&\displaystyle \quad 2\mathcal Q^{LL}(f_{L,0}^{n+1}, f_{L,1}^{n+1})
+ \underbrace{\mathcal Q_0^{LH}(f_{L,0}^n, f_{H,1}^n)}_{=O(\varepsilon^2 + \Delta t)} + \mathcal Q_0^{LH}(f_{L,1}^{n+1}, f_{H,0}^{n+1}) + 
\mathcal Q_1^{LH}(f_{L,0}^{n+1}, f_{H,0}^{n+1}) \notag\\[2pt]
&\displaystyle\label{SI_eps2} = v^L\cdot\nabla_x f_{L,0}^{n+1} + F_{L}^{n+1}\cdot\nabla_{v^L}f_{L,0}^{n+1} + O(\varepsilon^2). 
\end{align}
\\[2pt]
(\ref{SI_eps2}) is an equation for $f_{L,1}^{n+1}$ and can be equivalently written in terms of 
$$\phi_{L}^{n+1}=f_{L,1}^{n+1}\, (M_{L,0}^{n+1})^{-1}$$ according to: 
$$\Gamma_{L,0}\, \phi_{L}^{n+1} = -(M_{L,0}^{n+1})^{-1}\left[ v^L\cdot\nabla_x M_{L,0}^{n+1} + F_{L}^{n+1}\cdot\nabla_{v^L} M_{L,0}^{n+1}
-\mathcal Q_1^{LH}(M_{L,0}^{n+1}, M_{H,0}^{n+1})\right] + O(\varepsilon + \Delta t), $$
where $\Gamma_{L,0}$ is the linearized operator given by 
$$\Gamma_{L,0}\, \phi_{L}^{n+1} = (M_{L,0}^{n+1})^{-1} \left[2\mathcal Q^{LL}(M_{L,0}^{n+1},\, M_{L,0}^{n+1}\, \phi_{L}^{n+1}) + 
\mathcal Q_0^{LH}(M_{L,0}^{n+1}\, \phi_{L}^{n+1},\, M_{H,0}^{n+1})\right]. $$ 
As proved in \cite{Degond2}, the unique solution in $\left(\text{ker}(\Gamma_{L,0})\right)^{\perp}$ is given by
\begin{align*}
&\displaystyle\phi_{L}^{n+1} = \frac{1}{n_{L,0}^{n+1}}\left(-\left(\nabla_x n_{L,0}^{n+1} -\frac{F_{L}^{n+1}\, n_{L,0}^{n+1}}{T_{L,0}^{n+1}}\right)
\Psi_1(|v^L|) - n_{L,0}^{n+1}\frac{\nabla_x T_{L,0}^{n+1}}{T_{L,0}^{n+1}}\Psi_2(|v^L|) +\frac{n_{L,0}^{n+1}\, u_{H,0}^{n+1}}{T_{L,0}^{n+1}}\right) \cdot v^L \\[4pt]
&\displaystyle\qquad\quad + O(\varepsilon + \Delta t), 
\end{align*}
thus \begin{equation}\label{SI_fL1}f_{L,1}^{n+1} = \underbrace{M_{L,0}^{n+1}\, \phi_{L}^{n+1}}_{:= f_{L,1}^{\ast, n+1}} + 
O(\varepsilon + \Delta t). \end{equation}
\\[2pt]

We multiply (\ref{SI_APL}) by $\varepsilon$ and add it to (\ref{SI_num1N}), then get
\begin{align}
&\displaystyle \quad\frac{f_{L,0}^{n+1}-f_{L,0}^n}{\Delta t}+ \varepsilon\, \frac{f_{L,1}^{n+1}-f_{L,1}^n}{\Delta t}
+ \left(v^L\cdot\nabla_x f_{L,1}^{n} + F_{L}^{n}\cdot\nabla_{v^L}f_{L,1}^{n}\right) + \varepsilon \left(v^L\cdot\nabla_x f_{L,0}^{n} + F_{L}^{n}\cdot\nabla_{v^L}f_{L,0}^{n}\right) \notag\\[6pt]
&\displaystyle = \frac{1}{\varepsilon^2}\Big[\mathcal Q^{LL}(f_{L,0}^{n+1}, f_{L,0}^{n+1})+
\mathcal Q_0^{LH}(f_{L,0}^{n}, f_{H,0}^{n})\Big] \notag\\[4pt]
&\displaystyle\quad + \frac{1}{\varepsilon}\Big[2\mathcal Q^{LL}(f_{L,0}^{n+1}, f_{L,1}^{n+1})
+\mathcal Q_0^{LH}(f_{L,0}^n, f_{H,1}^n) + \mathcal Q_0^{LH}(f_{L,1}^{n+1}, f_{H,0}^{n+1}) + \mathcal Q_1^{LH}(f_{L,0}^{n+1}, f_{H,0}^{n+1}) \notag\\[4pt]
&\displaystyle\qquad\quad - \big(v^L\cdot\nabla_x f_{L,0}^{n+1} + F_{L}^{n+1}\cdot\nabla_{v^L}f_{L,0}^{n+1}\big)\Big]\notag\\[4pt]
&\displaystyle\quad + \mathcal Q^{LL}(f_{L,1}^n, f_{L,1}^n) +\mathcal Q_0^{LH}(f_{L,1}^n, f_{H,1}^n)
+ \mathcal Q_1^{LH}(f_{L,0}^n, f_{H,1}^n) + \mathcal Q_1^{LH}(f_{L,1}^n, f_{H,0}^n) + \mathcal Q_2^{LH}(f_{L,0}^{n+1}, f_{H,0}^{n+1})\notag\\[4pt]
&\displaystyle\label{SI_APL1}\quad + \varepsilon\Big[\mathcal Q_1^{LH}(f_{L,1}^n, f_{H,1}^n) + \mathcal Q_2^{LH}(f_{L,0}^n, f_{H,1}^n) 
+ \mathcal Q_2^{LH}(f_{L,1}^n, f_{H,0}^n) + \big(v^L\cdot\nabla_x f_{L,0}^{n+1} + F_{L}^{n+1}\cdot\nabla_{v^L}f_{L,0}^{n+1}\big)\Big]. 
\end{align}
Plugging in the leading order term of (\ref{SI_fL0}), (\ref{SI_fL1}) and comparing the $O(1)$ terms on both sides gives 
\begin{align}
&\displaystyle\quad\frac{M_{L,0}^{n+1} - M_{L,0}^n}{\Delta t} + v^L\cdot\nabla_x f_{L,1}^{\ast, n} + F_{L}^{n}\cdot\nabla_{v^L}f_{L,1^{\ast, n}}
\notag\\[4pt]
&\displaystyle 
= \mathcal Q^{LL}(f_{L,1}^{\ast, n}, f_{L,1}^{\ast, n}) + \mathcal Q_0^{LH}(f_{L,1}^{\ast, n}, f_{H,1}^{\ast, n})
+ \mathcal Q_1^{LH}(f_{L,0}^{\ast, n}, f_{H,1}^{\ast, n}) 
+ \mathcal Q_1^{LH}(f_{L,1}^{\ast, n}, M_{H,0}^{n}) \notag\\[4pt]
&\displaystyle\label{SI_APL1}\quad + \mathcal Q_2^{LH}(M_{L,0}^{n+1}, M_{H,0}^{n+1}) + O(\Delta t)\,. 
\end{align}
Integrate both sides of (\ref{SI_APL1}) against $1$, $v^L$, $|v^L|^2$ on $v^L$, 
by the statement 2(i) and (\ref{Q0_LH}) in Theorem \ref{property}, thus 
\begin{align*}
&\displaystyle \int \mathcal Q^{LL}(f_{L,1}^{\ast, n}, f_{L,1}^{\ast, n})\, |v^L|^2\, dv^L = \int \mathcal Q_0^{LH}(f_{L,1}^{\ast, n}, f_{H,1}^{\ast, n})\, |v^L|^2\, dv^L = 0, \\[4pt]
&\displaystyle  \int  \mathcal Q_1^{LH}(M_{L,0}^{n}, f_{H,1}^{\ast, n})\, |v^L|^2\, dv^L = \int \mathcal Q_0^{HL}(f_{H,1}^{\ast, n}, M_{L,0}^{n})\, |v^L|^2\, dv^L = 0\,. 
\end{align*}
\\[2pt]
Integrals of $\frac{M_{L,0}^{n+1} - M_{L,0}^n}{\Delta t}$ are 
$$\frac{d}{dt}\left(n_{L,0}^{n}, \, n_{L,0}^{n} u_{L,0}^{n}, \, n_{L,0}^{n}(\frac{1}{2}|u_{L,0}^{n}|^2 +\frac{3}{2}T_{L,0}^{n})\right)^{T}. $$
Analogous to the calculation in \cite{Degond2}, then 
\begin{align*}
&\displaystyle\quad \int\left[\mathcal Q_1^{LH}(f_{L,1}^{\ast, n}, M_{H,0}^{n}) +  \mathcal Q_2^{LH}(M_{L,0}^{n+1}, M_{H,0}^{n+1})\right] |v^L|^2\, dv^L \\[4pt]
&\displaystyle = \int \left[\mathcal Q_0^{HL}(M_{H,0}^{n}, f_{L,1}^{\ast, n}) + \mathcal Q_1^{HL}(M_{H,0}^{n+1}, M_{L,0}^{n+1})\right] |v^H|^2\, dv^H \\[4pt]
&\displaystyle = \int \left[\mathcal Q_0^{HL}(M_{H,0}^{n+1}, f_{L,1}^{\ast, n+1}) + \mathcal Q_1^{HL}(M_{H,0}^{n+1}, M_{L,0}^{n+1})\right] |v^H|^2\, dv^H + O(\Delta t) \\[4pt]
&\displaystyle = u_{H,0}^{n+1}\cdot \left[\nabla_x(n_{L,0}^{n+1}\, T_{L,0}^{n+1}) -F_{L}^{n+1}\, n_{L,0}^{n+1}\right] +
3\, \frac{\lambda(T_{L,0}^{n+1})}{T_{L,0}^{n+1}}\, n_{L,0}^{n+1}\, n_{H,0}^{n+1}\, (T_{H,0}^{n+1}- T_{L,0}^{n+1}) + O(\Delta t). 
\end{align*}
Therefore, the limit of our scheme is given by 
\begin{align}
&\displaystyle\frac{\partial n_{L,0}^n}{\partial t} + \nabla_x\cdot (n_{L,0}^n\, u_{H,0}^n) \notag\\[4pt]
&\displaystyle\label{SI_Limit1}\qquad\quad -\nabla_x\cdot\left[D_{11}\left(\nabla_x n_{L,0}^n -\frac{F_{L}^n\, n_{L,0}^n}{T_{L,0}^n}\right) + 
D_{12}\left(n_{L,0}^n\, \frac{\nabla_x T_{L,0}^n}{T_{L,0}^n}\right)\right]= O(\Delta t),   \\[6pt]
&\displaystyle\frac{\partial}{\partial t}\left(\frac{3}{2}n_{L,0}^n\, T_{L,0}^n\right) + \nabla_x\cdot \left(\frac{5}{2}n_{L,0}^n\, T_{L,0}^n\, u_{H,0}^n\right) 
- n_{L,0}^n\, F_{L}^n\, u_{H,0}^n \notag\\[4pt]
&\displaystyle\qquad\qquad \qquad\qquad -\nabla_x \cdot\left[D_{21}\left(\nabla_x n_{L,0}^n-\frac{F_L^n\, n_{L,0}^n}{T_{L,0}^n}\right) + 
D_{22}\left(n_{L,0}^n\, \frac{\nabla_x T_{L,0}^n}{T_{L,0}^n}\right)\right] \notag\\[4pt]
&\displaystyle\qquad\qquad \qquad\qquad  + F_L^{n}\cdot\left[D_{11}\left(\nabla_x n_{L,0}^n-\frac{F_L^n\, n_{L,0}^n}{T_{L,0}^n}\right) + D_{12}
\left(n_{L,0}^n\, \frac{\nabla_x T_{L,0}^n}{T_{L,0}^n}\right)\right]\notag\\[4pt]
&\displaystyle\label{SI_Limit2}\quad = u_{H,0}^{n+1}\cdot \left[\nabla_x(n_{L,0}^{n+1}\, T_{L,0}^{n+1}) -F_{L}^{n+1}\, n_{L,0}^{n+1}\right] +
3\, \frac{\lambda(T_{L,0}^{n+1})}{T_{L,0}^{n+1}}\, n_{L,0}^{n+1}\, n_{H,0}^{n+1}\, (T_{H,0}^{n+1}- T_{L,0}^{n+1}) +  O(\Delta t). 
\end{align}
This is first order (in $\Delta t$) consistent to the the implicit numerical discretization of the limit equation (\ref{nL0})--(\ref{nL1}). 
\\[2pt]

Next we look at the system for the heavy particles. Inserting the expansion $$\mathcal Q_0^{HL} +\varepsilon \mathcal Q_1^{HL}+ \varepsilon^2 \mathcal Q_2^{HL} + O(\varepsilon^3) $$
into (\ref{SI_num4N}), one has
\begin{align}
&\displaystyle\quad \frac{f_{H,1}^{n+1}-f_{H,1}^n}{\Delta t}  
 + \big(v^H\cdot\nabla_x f_{H,0}^{n} + F_{H}^{n}\cdot\nabla_{v^H}f_{H,0}^{n}\big)\notag\\[4pt]
&\displaystyle = \frac{1}{\varepsilon}\Big[\mathcal Q_0^{HL}(f_{H,0}^{n+1}, f_{L,1}^{n+1}) + \mathcal Q_0^{HL}(f_{H,1}^n, f_{L,0}^n) 
+ 2\mathcal Q^{HH}(f_{H,0}^{n+1}, f_{H,1}^{n+1}) + \mathcal Q_1^{HL}(f_{H,0}^{n+1}, f_{L,0}^{n+1})\notag\\[4pt]
&\displaystyle\qquad - \big(v^H\cdot\nabla_x f_{H,0}^{n+1} + F_{H}^{n}\cdot\nabla_{v^H}f_{H,0}^{n+1}\big)\Big] \notag\\[4pt]
&\displaystyle\quad + \mathcal Q^{HH}(f_{H,1}^n, f_{H,1}^n) +\mathcal Q_1^{HL}(f_{H,0}^{n+1}, f_{L,1}^{n+1})
+\mathcal Q_1^{HL}(f_{H,1}^n, f_{L,0}^n) + \mathcal Q_0^{HL}(f_{H,1}^n, f_{L,1}^n) + \mathcal Q_2^{HL}(f_{H,0}^{n+1}, f_{L,0}^{n+1})\notag\\[4pt]
&\displaystyle\label{AP_H2}\quad + \big(v^H\cdot\nabla_x f_{H,0}^{n+1} + F_{H}^{n}\cdot\nabla_{v^H}f_{H,0}^{n+1}\big) 
+ \varepsilon\Big[\mathcal Q_1^{HL}(f_{H,1}^n, f_{L,1}^n) + \mathcal Q_2^{HL}(f_{H,0}^{n+1}, f_{L,1}^{n+1})+\mathcal Q_2^{HL}(f_{H,1}^{n+1}, f_{L,0}^{n+1})\Big]. 
\end{align}
\\[2pt]
We multiply (\ref{AP_H2}) by $\varepsilon$ and add it up with (\ref{SI_num3N}), then get 
\begin{align}
&\displaystyle\quad \frac{f_{H,0}^{n+1}-f_{H,0}^n}{\Delta t} + \varepsilon\, \frac{f_{H,1}^{n+1}-f_{H,1}^n}{\Delta t}  
+ \varepsilon \big(v^H\cdot\nabla_x f_{H,1}^{n} + F_{H}^{n}\cdot\nabla_{v^H}f_{H,1}^{n}\big) + 
\big(v^H\cdot\nabla_x f_{H,0}^{n} + F_{H}^{n}\cdot\nabla_{v^H}f_{H,0}^{n}\big) \notag\\[4pt]
&\displaystyle = \frac{1}{\varepsilon}\left[\mathcal Q^{HH}(f_{H,0}^{n+1}, f_{H,0}^{n+1})+
\mathcal Q_0^{HL}(f_{H,0}^{n}, f_{L,0}^{n})\right] \notag\\[4pt]
&\displaystyle\quad + \mathcal Q_0^{HL}(f_{H,0}^{n+1}, f_{L,1}^{n+1}) + \mathcal Q_0^{HL}(f_{H,1}^n, f_{L,0}^n) 
+ 2\mathcal Q^{HH}(f_{H,0}^{n+1}, f_{H,1}^{n+1}) + \mathcal Q_1^{HL}(f_{H,0}^{n+1}, f_{L,0}^{n+1}) \notag\\[4pt]
&\displaystyle\quad - \big(v^H\cdot\nabla_x f_{H,0}^{n+1} + F_{H}^{n}\cdot\nabla_{v^H}f_{H,0}^{n+1}\big) \notag\\[4pt]
&\displaystyle\quad + \varepsilon\Big[\mathcal Q^{HH}(f_{H,1}^n, f_{H,1}^n) +\mathcal Q_1^{HL}(f_{H,0}^{n+1}, f_{L,1}^{n+1})
+\mathcal Q_1^{HL}(f_{H,1}^n, f_{L,0}^n) + \mathcal Q_0^{HL}(f_{H,1}^n, f_{L,1}^n) + \mathcal Q_2^{HL}(f_{H,0}^{n+1}, f_{L,0}^{n+1})\notag\\[4pt]
&\displaystyle\qquad\quad + \big(v^H\cdot\nabla_x f_{H,0}^{n+1} + F_{H}^{n}\cdot\nabla_{v^H}f_{H,0}^{n+1}\big)\Big] \notag\\[4pt]
&\displaystyle\label{AP_H}\quad + \varepsilon^2 \left[\mathcal Q_1^{HL}(f_{H,1}^n, f_{L,1}^n) + \mathcal Q_2^{HL}(f_{H,0}^{n+1}, f_{L,1}^{n+1}) +\mathcal Q_2^{HL}(f_{H,1}^{n+1}, f_{L,0}^{n+1})\right]. 
\end{align}
Plugging in the leading order term of (\ref{SI_fH0}) and comparing the $O(1)$ terms on both sides, one gets 
\begin{align}
&\displaystyle\frac{M_{H,0}^{n+1}- M_{H,0}^n}{\Delta t} = 
2\mathcal Q^{HH}(M_{H,0}^{n+1}, f_{H,1}^{\ast, n+1}) + \mathcal Q_0^{HL}(M_{H,0}^{n+1}, f_{L,1}^{\ast, n+1}) + 
\underbrace{\mathcal Q_0^{HL}(f_{H,1}^{\ast, n}, M_{L,0}^{n})}_{=0}\notag\\[2pt]
&\displaystyle\label{SI_eps0}\qquad\qquad\qquad\quad + \mathcal Q_1^{HL}(M_{H,0}^{n+1}, M_{L,0}^{n+1}) 
- \left(v^H\cdot\nabla_x M_{H,0}^{n} + F_{H}^{n}\cdot\nabla_{v^H}M_{H,0}^{n}\right) + O(\Delta t). 
\end{align}
(\ref{SI_eps0}) can be equivalently written for $$\phi_{H,1}^{n+1}=(M_{H,0}^{n+1})\, f_{H,1}^{\ast, n+1}$$ with 
\begin{equation}
\label{SI_Gamma}\Gamma_{H,0}\, \phi_{H,1}^{n+1} = (M_{H,0}^{n+1})^{-1}\, S_{H,1}^{n+1} + O(\Delta t). \end{equation}
$\Gamma_{H,0}$ is a linearization operator given by 
$$\Gamma_{H,0}\, \phi_{H,1}^{n+1} = 2\, (M_{H,0}^{n+1})^{-1}\, \mathcal Q^{HH}(M_{H,0}^{n+1}, \, M_{H,0}^{n+1}\, \phi_{H,1}^{n+1}), $$
and $S_{H,1}^{n+1}$ is 
\begin{align}
&\displaystyle S_{H,1}^{n+1} = \left(D_t + v^H\cdot\nabla_x + F_{H}^{n}\cdot\nabla_{v^H}\right) M_{H,0}^{n} \notag\\[4pt]
&\displaystyle\qquad\qquad - \mathcal Q_0^{HL}(M_{H,0}^{n+1}, \, f_{L,1}^{n+1}) - \mathcal Q_1^{HL}(M_{H,0}^{n+1}, M_{L,0}^{n+1}). 
\end{align}
The necessary and sufficient condition of solvability of equation (\ref{SI_Gamma}) is
\begin{equation} \label{SI_Cond}\int_{\mathbb R^3}\, S_{H,1}^{n+1}\begin{pmatrix} 1\\ v^H \\|v^H|^2\end{pmatrix} dv^H = 
O(\Delta t)\, \mathbb I_{3}\,. 
\end{equation}

The following is analogous to the proof shown in \cite{Degond2}, except that we have a discrete counterpart here. 
With details omitted, (\ref{SI_Cond}) thus gives 
\begin{align}
&\displaystyle\label{Limit_H1} \frac{\partial n_{H,0}^n}{\partial t} + \nabla_x\cdot (n_{H,0}^n\, u_{H,0}^n) = O(\Delta t),    \\[8pt]
&\displaystyle \frac{\partial}{\partial t}(n_{H,0}^n\, u_{H,0}^n)  + \nabla_x\cdot(n_{H,0}^n\, u_{H,0}^n \otimes u_{H,0}^n) + 
\nabla_x (n_{H,0}^n\, T_{H,0}^n) - n_{H,0}^n\, F_{H}^n \notag\\[4pt]
&\displaystyle\label{Limit_H2}\quad = - \big(\nabla_x (n_{L,0}^{n+1}\, T_{L,0}^{n+1}) - F_{L}^{n+1}\, n_{L,0}^{n+1}\big) + O(\Delta t),  \\[8pt]
&\displaystyle \frac{\partial}{\partial t}\left(\frac{n_{H,0}^n\, |u_{H,0}^n|^2}{2} + \frac{3}{2}n_{H,0}^n\, T_{H,0}^n\right)
+\nabla_x \cdot\left(\left(\frac{n_{H,0}^n\, |u_{H,0}^n|^2}{2} +\frac{5}{2} n_{H,0}^n\, T_{H,0}^n\right)u_{H,0}^n\right) 
- n_{H,0}^n\, F_{H}^n\cdot u_{H,0}^n \notag\\[4pt]
&\displaystyle \quad  = -u_{H,0}^{n+1} \cdot\left[\nabla_x (n_{L,0}^{n+1}\, T_{L,0}^{n+1}) - F_{L}^{n+1}\, n_{L,0}^{n+1}\right] 
- 3\, \frac{\lambda(T_{L,0}^{n+1})}{T_{L,0}^{n+1}}\, n_{L,0}^{n+1}\, n_{H,0}^{n+1}\, (T_{H,0}^{n+1} - T_{L,0}^{n+1}) \notag\\[4pt]
&\displaystyle\label{Limit_H3} \qquad + O(\Delta t). 
\end{align}
Therefore, (\ref{Limit_H1}), (\ref{Limit_H2}) and (\ref{Limit_H3}) are consistent with the discrete scheme of the hydrodynamic limit system 
(\ref{nH0})--(\ref{nH2}), up to a numerical error of $O(\Delta t)$. 
\subsection{Splitting of convection from the collision}
As in \cite{JPT1, JPT2}, we adopt a first-order time splitting approach to separate the convection from the collision operators.
To summarize, our scheme is given by the following equations: 
\\[8pt]
{\bf Moment equations for $f_{L,0}$ and $f_{H,0}$}: 
\begin{align}
&\displaystyle\label{SI_L0_a} (P_0)_{L,0}^{n+1} = (P_0)_{L,0}^n + \Delta t \int_{\mathbb R^3} v^L\cdot\nabla_{x}f_{L,1}^{n}\, dv^L, 
 \\[2pt]
&\displaystyle\label{SI_L0_b} (P_1)_{L,0}^{n+1} = (P_1)_{L,0}^n + \frac{\Delta t}{\varepsilon^2}\int_{\mathbb R^3} \phi_1^L\, 
 \mathcal Q_0^{LH}(f_{L,0}^n, f_{H,0}^n)(v^L)\, dv^L + \Delta t \int_{\mathbb R^3} \phi_1^L\, v^L\cdot\nabla_{x}f_{L,1}^{n}\, dv^L, \\[2pt]
&\displaystyle\label{SI_L0_c} (P_2)_{L,0}^{n+1} = (P_2)_{L,0}^n + \Delta t \int_{\mathbb R^3} \phi_2^L\, v^L\cdot\nabla_{x}f_{L,1}^{n} \, dv^L, 
\\[4pt]
&\displaystyle\label{SI_H0_a} (P_0)_{H,0}^{n+1} = (P_0)_{H,0}^n + \varepsilon\Delta t \int_{\mathbb R^3} v^H \cdot\nabla_{x}f_{H,1}^{n}\, dv^H, 
 \\[2pt]
&\displaystyle\label{SI_H0_b} (P_i)_{H,0}^{n+1} = (P_i)_{H,0}^n + \frac{\Delta t}{\varepsilon^2}\int_{\mathbb R^3}
\phi_i^H\, \mathcal Q_0^{HL}(f_{H,0}^n, f_{L,0}^n)(v^H)\, dv^H +  \varepsilon\Delta t \int_{\mathbb R^3}\phi_i^H\, v^H \cdot\nabla_{x}f_{H,1}^{n}\, dv^H, 
\end{align}
where $\phi_i^L$, $\phi_i^H$ are defined in (\ref{Phi}) and $i=1, 2$. 
\\[10pt]
\indent The scheme for $f_{L,0}$, $f_{L,1}$, $f_{H,0}$, $f_{H,1}$ are given by:
\\[4pt]
\noindent{\bf Step 1: The implicit collision step}
\begin{align}
&\displaystyle \label{SI_NL0} \frac{f_{L,0}^{\ast}-f_{L,0}^{n}}{\Delta t} =\frac{1}{\varepsilon^2}\bigg[\mathcal Q^{LL}(f_{L,0}^{n}, f_{L,0}^{n}) - \mathcal P(f_{L,0}^n) + \mathcal P(f_{L,0}^{\ast})
+\mathcal Q_0^{LH}(f_{L,0}^{n}, f_{H,0}^{n})\bigg], \\[6pt]
&\displaystyle\frac{f_{L,1}^{\ast}-f_{L,1}^n}{\Delta t} = \frac{1}{\varepsilon^2}\bigg[\frac{1}{\varepsilon}\left(\mathcal Q^{LH}_{\varepsilon}(f_{L,0}^{\ast}, f_{H,0}^{\ast}) -\mathcal Q_0^{LH}(f_{L,0}^{\ast}, f_{H,0}^{\ast})\right) \notag\\[2pt]
&\displaystyle\qquad\qquad\qquad\quad + \frac{1}{2}\bigg[\mathcal Q^{LL}(f_{L,0}^{n}+f_{L,1}^{n}, f_{L,0}^{n}+f_{L,1}^{n}) + \mu (f_{L,0}^n+f_{L,1}^n) - \mu (f_{L,0}^{\ast}+f_{L,1}^{\ast}) \notag\\[2pt]
&\displaystyle \qquad\qquad\qquad\qquad - \left(\mathcal Q^{LL}(f_{L,0}^{n}-f_{L,1}^{n}, f_{L,0}^{n}-f_{L,1}^{n}) + \mu (f_{L,0}^n-f_{L,1}^n) - \mu (f_{L,0}^{\ast}-f_{L,1}^{\ast})\right)\bigg] \notag\\[2pt]
&\displaystyle\qquad\qquad\qquad\quad  + \varepsilon\mathcal Q^{LL}(f_{L,1}^{n}, f_{L,1}^{n})  + \mathcal Q^{LH}_{\varepsilon}(f_{L,0}^{n}, f_{H,1}^{n}) + \left(\mathcal Q^{LH}_{\varepsilon}(f_{L,1}^{n}, f_{H,0}^{n}) - \mathcal Q_0^{LH}(f_{L,1}^{n}, f_{H,0}^{n})\right) \notag\\[2pt]
&\displaystyle\label{SI_NL1}\qquad\qquad\qquad\quad + \mathcal Q_0^{LH}(f_{L,1}^{\ast}, f_{H,0}^{\ast}) +\varepsilon\mathcal Q^{LH}_{\varepsilon}(f_{L,1}^{n}, f_{H,1}^n) - (1-\varepsilon^2\psi_1)
\left(v^L\cdot\nabla_x f_{L,0}^{\ast} + F^L\cdot\nabla_{v^L}f_{L,0}^{\ast}\right)\bigg], 
\end{align}
\begin{align}
&\displaystyle\label{SI_NH0}\frac{f_{H,0}^{\ast}-f_{H,0}^n}{\Delta t} = \frac{1}{\varepsilon}\bigg[\mathcal Q^{HH}(f_{H,0}^{\ast}, f_{H,0}^{\ast}) - \mathcal P(f_{H,0}^n) + \mathcal P(f_{H,0}^{\ast}) + \mathcal Q_0^{HL}(f_{H,0}^{n}, f_{L,0}^{n})\bigg], \\[6pt]
&\displaystyle\frac{f_{H,1}^{\ast}-f_{H,1}^n}{\Delta t} = \frac{1}{\varepsilon}\bigg[\frac{1}{\varepsilon}\left(\mathcal Q^{HL}_{\varepsilon}(f_{H,0}^{\ast}, f_{L,0}^{\ast})-\mathcal Q_0^{HL}(f_{H,0}^{\ast}, f_{L,0}^{\ast})\right)\notag\\[2pt]
&\displaystyle\qquad\qquad\qquad\quad
+ \frac{1}{2}\bigg[\mathcal Q^{HH}(f_{H,0}^{n}+f_{H,1}^{n}, f_{H,0}^{n}+f_{H,1}^{n}) + \mu(f_{H,0}^n+f_{H,1}^n) 
- \mu(f_{H,0}^{\ast}+f_{H,1}^{\ast}) \notag\\[2pt]
&\displaystyle\qquad\qquad\qquad\qquad - \left(\mathcal Q^{HH}(f_{H,0}^{n}-f_{H,1}^{n}, f_{H,0}^{n}-f_{H,1}^{n}) + \mu(f_{H,0}^n-f_{H,1}^n)
- \mu(f_{H,0}^{\ast}-f_{H,1}^{\ast})\right)\bigg] \notag \\[2pt]
&\displaystyle\qquad\qquad\qquad\quad  + \varepsilon\mathcal Q^{HH}(f_{H,1}^{n}, f_{H,1}^{n}) + \mathcal Q^{HL}_{\varepsilon}(f_{H,0}^{\ast}, f_{L,1}^{\ast}) + \mathcal Q^{HL}_{\varepsilon}(f_{H,1}^{n}, f_{L,0}^{n}) +\varepsilon\mathcal Q^{HL}_{\varepsilon}(f_{H,1}^{n}, f_{L,1}^n)\notag\\[2pt]
&\displaystyle\label{SI_NH1} \qquad\qquad\qquad\quad {\color{black} - (1-\varepsilon\psi_2)\Big(v^H\cdot\nabla_x f_{H,0}^{\ast} + 
F_{H}^{n}\cdot\nabla_{v^H}f_{H,0}^{\ast}\Big)\bigg].}
\end{align}
The order is to first solve (\ref{SI_NL0}), (\ref{SI_NH0}), then solve (\ref{SI_NL1}) and (\ref{SI_NH1}). 
\\[12pt]
{\bf Step 2: The explicit transport step}
\begin{align}
&\displaystyle \frac{f_{L,0}^{n+1}-f_{L,0}^{\ast}}{\Delta t}+\left(v^L\cdot\nabla_x f_{L,1}^{\ast} + F_{L}^{\ast}\cdot\nabla_{v^L}f_{L,1}^{\ast}\right) =0\,, \\[4pt]
&\displaystyle \frac{f_{L,1}^{n+1}-f_{L,1}^{\ast}}{\Delta t}+ {\color{black}\psi_1}\left(v^L\cdot\nabla_x f_{L,0}^{\ast} + F_{L}^{\ast}\cdot\nabla_{v^L}f_{L,0}^{\ast}\right)=0 \,.
\end{align}
and 
\begin{align}
&\displaystyle \frac{f_{H,0}^{n+1}-f_{H,0}^{\ast}}{\Delta t}+ \varepsilon\left(v^H\cdot\nabla_x f_{H,1}^{\ast} + F_{H}^{\ast}\cdot\nabla_{v^H}f_{H,1}^{\ast}\right) =0\,, \\[4pt]
&\displaystyle \frac{f_{H,1}^{n+1}-f_{H,1}^{\ast}}{\Delta t}+ {\color{black}\psi_2}\left(v^H\cdot\nabla_x f_{H,0}^{\ast} + F_{H}^{\ast}\cdot\nabla_{v^H}f_{H,0}^{\ast}\right)=0\,, 
\end{align}
where {\color{black}$$\psi_1=\min\{1, \frac{1}{\varepsilon^2}\}, \qquad \psi_2=\min\{1, \frac{1}{\varepsilon}\}. $$}

\section{Conclusion and future work}
\label{sec4}

In this paper, we develop asymptotic-preserving time discretizations 
for disparate mass binary gas or  plasma for both the homogeneous and inhomogeneous cases, 
at the relaxation time scale, for both the Boltzmann and Fokker-Planck-Landau collision operators.
We introduce a novel splitting of the system and a carefully designed
explicit-implicit time discretization so to first guarantee the correct
asymptotic behavior at the relaxation time limit and also significantly reduces
the algebraic complexity which will be comparable to their single species counterparts.
The design of the AP schemes are strongly guided by the asymptotic
behavior of the system studied in \cite{Degond, Degond2}. We also prove
that a simplied version of the time discretization is asymptotic-preserving.

In the follow-up work, spatial and velocity discretizations will be discussed,
along with extensive numerical simulations and experiments.
Moreover, we plan to address the issue of uncertainty quantification (UQ),
by adding random inputs into the system, and develop efficient numerical methods for such uncertain kinetic system. 

\section*{Acknowledgement}
The third author would like to thank Dr. Ruiwen Shu for a helpful discussion. We thank the referees for their helpful comments. 
\\[10pt]

\chapter{{\bf\Large{ Appendix}}}
\renewcommand{\theequation}{A.\arabic{equation}}
\\[2pt]

{\bf Definitions of $\mathcal Q$}
\\[6pt]
In the Fokker-Planck-Landau case, the collision operators 
 are given by 
\begin{align*}
&\displaystyle   \mathcal Q_{\mathcal {L}}^{LL} :=\mathcal Q_{\mathcal {L}}^{LL}(f^L, f^L)(v^L)=\nabla_{v^L}\cdot\int_{\mathbb R^3}\, B(v^L-v_{\ast}^L)\, S(v^L-v_{\ast}^L)(\nabla_{v^L}f^L f_{\ast}^L-\nabla_{v_{\ast}^L}f^L f^L)\, dv_{\ast}^L, \\[2pt]
&\displaystyle  \mathcal Q_{\mathcal {L}}^{HH}:=\mathcal Q_{\mathcal {L}}^{HH}(f^H, f^H)(v^H)=\nabla_{v^H}\cdot\int_{\mathbb R^3}\, B(v^H-v_{\ast}^H)\, S(v^H-v_{\ast}^H)(\nabla_{v^H}f^H f_{\ast}^H-\nabla_{v_{\ast}^H}f^H f^H)\, dv_{\ast}^H, \\[2pt]
&\displaystyle  \mathcal Q_{\mathcal {L},\varepsilon}^{LH}:=\mathcal Q_{\mathcal {L},\varepsilon}^{LH}(f^L, f^H(\varepsilon v^H))(v^L)=(1+\varepsilon^2)^{\frac{\gamma+2}{2}}\, \nabla_{v^L}\cdot\int_{\mathbb R^3}\, 
B(\frac{v^L-\varepsilon v^H}{\sqrt{1+\varepsilon^2}})\, S(v^L-\varepsilon v^H)(\nabla_{v^L}f^L f^H -\varepsilon\nabla_{v^H}f^H f^L)\, dv^H, 
\\[2pt]
&\displaystyle \mathcal Q_{\mathcal {L},\varepsilon}^{HL}:=\mathcal Q_{\mathcal {L},\varepsilon}^{HL}(f^H(\varepsilon v^H), f^L)(v^H)=(1+\varepsilon^2)^{\frac{\gamma+2}{2}}\, \nabla_{v^H}\cdot\int_{\mathbb R^3}\, 
B(\frac{v^L-\varepsilon v^H}{\sqrt{1+\varepsilon^2}})\, S(v^L-\varepsilon v^H)(\nabla_{v^L}f^L f^H -\varepsilon\nabla_{v^H}f^H f^L)\, dv^L,  
\end{align*}
where the matrix $S(w)$ and the intra-molecular potential $B(w)$, respectively, are given by  
$$S(w)= \text{Id} - \frac{w\otimes w}{|w|^2}, \qquad B(w)=\frac{1}{2} |w|^{\gamma+2}. $$
In particular, $B\left(\frac{1}{\sqrt{1+\varepsilon^2}}w \right)= \frac{1}{2} (1+\varepsilon^2)^{-\frac{\gamma+2}{2}}|w|^{\gamma+2}$, and 
the value  $\gamma=-3$ corresponds to Coulomb interactions.   \\ 
\\[2pt]

In the Boltzmann case, the collision operators in center of mass -- relative velocity coordinates expressed in the angular scattering direction 
$\sigma$, are given by
\begin{align*}
&\displaystyle \mathcal Q_{\mathcal {B}}^{LL}:= \mathcal Q^{LL}(f^L, f^L)(v^L) = \int_{\mathbb R^3}\int_{S^2}\, B^L(v^L-v_{\ast}^L,\sigma)\, (f^{\prime, L}f_{\ast}^{\prime, L} - f^L f_{\ast}^L)\, d\sigma dv_{\ast}^L, \\[2pt]
&\displaystyle \mathcal Q_{\mathcal {B}}^{HH}:=\mathcal Q^{HH}(f^H, f^H)(v^H) = \int_{\mathbb R^3}\int_{S^2}\, B^H(v^H-v_{\ast}^H, \sigma)\, (f^{\prime, H}f_{\ast}^{\prime, H} - f^H f_{\ast}^H)\, d\sigma dv_{\ast}^H,  \\[2pt]
&\displaystyle \mathcal Q^{LH}_{\mathcal {B},\varepsilon} := \mathcal Q^{LH}_{\varepsilon}(f^L, f^H(\varepsilon v^H))(v^L)=(1+\varepsilon^2)^{\frac{\gamma}{2}}\, \int_{\mathbb R^3}\int_{S^2}\, B(\frac{v^L-\varepsilon v^H}{\sqrt{1+\varepsilon^2}}, \sigma)\,
(f^{\prime L, \varepsilon}f^{\prime H, \varepsilon}-f^L f^H)\, d\sigma dv^H, \\[2pt]
&\displaystyle \mathcal Q^{HL}_{\mathcal {B}, \varepsilon} := \mathcal Q^{HL}_{\varepsilon}(f^H(\varepsilon v^H), f^L)(v^H)=
\left(\frac{1+\varepsilon^2}{\varepsilon^2}\right)^{\frac{\gamma}{2}} \int_{\mathbb R^3}\int_{S^2}\,  B(\frac{\varepsilon}{\sqrt{1+\varepsilon^2}}(v^L-\varepsilon v^H), \sigma)\,
(f^{\prime L, \varepsilon}f^{\prime H, \varepsilon}-f^L f^H)\, d\sigma dv^L, 
\end{align*}
with $$v^{\prime L, \varepsilon}=v^L + \frac{1}{1+\varepsilon^2}\left(|v^L - \varepsilon v^H| \sigma - (v^L - \varepsilon v^H) \right)
= \frac{\varepsilon^2 v^L + \varepsilon v^H + |v^L - \varepsilon v^H| \sigma}{1+\varepsilon^2}, 
$$
and $$ v^{\prime H, \varepsilon}=\varepsilon v^H - \frac{\varepsilon^2}{1+\varepsilon^2}\left(|v^L - \varepsilon v^H| \sigma - (v^L - \varepsilon v^H) \right)
= \frac{\varepsilon^2 v^L + \varepsilon v^H - \varepsilon^2 |v^L - \varepsilon v^H| \sigma}{1 + \varepsilon^2}. 
$$
Here the collision kernel $B$ is assumed to be in the form $$ B(w,\sigma)  = \frac12 |w|^{\gamma}\, b(\frac{w}{|w|}\cdot\sigma). $$  
\\[2pt]

{\bf The penalty methods}
\\[6pt]
For the Boltzmann equation, the best choice of this relaxation operator shown in \cite{Filbet-Jin} is
$$P(f)=\beta\left(\mathcal M_{\rho, u, T} - f\right), $$
where $\beta>0$ is an upper bound of $||\nabla\mathcal Q(\mathcal M_{\rho, u, T})||$. 
Another simple example of $\beta$ at time $t^n$ is 
$$\beta^n =\sup \left|\frac{\mathcal Q(f^n, f^n)-\mathcal Q(f^{n-1}, f^{n-1})}{f^n - f^{n-1}}\right|. $$

We briefly review the penalty method introduced in \cite{Filbet-Jin} for the Boltzmann equation in the form: 
$$ \partial_t f + v\cdot \nabla_x f = \frac{1}{\varepsilon}\mathcal Q_{B}(f, f), $$
the discretized scheme is given by 
\begin{equation}\label{FJ} \frac{f^{n+1}-f^n}{\Delta t} + v \cdot\nabla_x f^n = \frac{\mathcal Q_{B}(f^n, f^n) - P(f^n)}{\varepsilon} + \frac{P(f^{n+1})}{\varepsilon}, 
\end{equation}
where $P(f)=\beta\left[ \mathcal M_{\rho, u, T}(v) - f(v)\right]$. 
Multiplying (\ref{FJ}) by $\phi(v)=(1, v, |v|^2)^{T}$, one gets the macroscopic 
quantities $U:=(\rho, \rho u, T)$: 
$$U^{n+1} = \int \phi(v)(f^n - \Delta t v\cdot\nabla_x f^n)\, dv. $$
$U^{n+1}$ is obtained explicitly, which defines $\mathcal M^{n+1}$, thus $f^{n+1}$ can be computed explicitly. 
\\[2pt]

On the other hand, \cite{JinYan} discusses the penalty method for solving the multiscale Fokker-Planck-Landau equation: 
\begin{equation}\label{toy}\partial_t f + v\cdot \nabla_x f = \frac{1}{\varepsilon}\mathcal Q_{L}(f, f),  \end{equation}
The authors in \cite{JinYan} demonstrate analytically and numerically that the best choice of the penalization operator
is the linear Fokker-Planck (FP) operator, 
\begin{equation}\label{FP}P_{FP}(f) = \nabla_v \cdot\left(M\nabla_v\left(\frac{f}{M}\right)\right), \end{equation}
where $$ M(x,v)=\frac{\rho(x)}{(2\pi T(x))^{N_v/2}}\exp\left(-\frac{(v-u(x))^2}{2T(x)}\right). $$
The first order AP scheme for (\ref{toy}) is given by 
$$\frac{f^{n+1}-f^n}{\Delta t}=\frac{1}{\varepsilon}\left(\mathcal Q(f^n, f^n) - \beta P^n f^n +\beta P^{n+1}f^{n+1}\right), $$
where $\beta$ is chosen large enough to ensure stability. For example, let $\displaystyle\beta=\beta_0 \max_{v}\lambda(D_A(f))$, 
with $\beta_0>\frac{1}{2}$ and $\lambda(D_A)$ is the spectral radius of the positive symmetric matrix $D_A$, defined by 
$$D_A(f) = \int_{\mathbb R^3} A(v-v_{\ast})f_{\ast}\, dv_{\ast}, $$
with $$A(z)=|z|^{\gamma+2}\left(\mathbb I -\frac{z\otimes z}{|z|^2}\right). $$

Compute the moments of $f^n$ by 
$$ (\rho, \rho u , \rho T)^{n+1} = \int_{\mathbb R^3} \left(1, v, \frac{(v-u)^2}{2}\right) f^n dv, $$
and update $M^{n+1}$. One can then solve $f^{n+1}$ by 
\begin{equation}
\label{f_FP1} f^{n+1} = \left(1 - \frac{\beta \Delta t}{\varepsilon}P^{n+1}\right)^{-1}
\left( f^n + \frac{\Delta t}{\varepsilon}(\mathcal Q(f^n) - \beta P^n f^n)\right). 
\end{equation}

Introduce the symmetrized operator \cite{JinYan}
$$ \tilde P h = \frac{1}{\sqrt{M}}\nabla_v \cdot\left(M\nabla_v \left(\frac{h}{\sqrt{M}}\right)\right), $$
then the penalty operator is $P_{FP} f = \sqrt{M}\tilde P \left(\frac{f}{\sqrt{M}}\right)$. 
Rewrite (\ref{f_FP1}) as
$$ \left(\frac{f}{\sqrt{M}}\right)^{n+1} =  \left(1 - \frac{\beta \Delta t}{\varepsilon}P^{n+1}\right)^{-1}\bigg\{\frac{1}{\sqrt{M^{n+1}}}
\left[f^n +\frac{\Delta t}{\varepsilon}\left(\mathcal Q(f^n) - \beta \sqrt{M^n}\tilde P^n \left(\frac{f^n}{\sqrt{M^n}}\right)\right)\right]\bigg\}. $$

The discretization of $\tilde P$ in one dimension is given by 
\begin{align*}
&\displaystyle (\tilde P h)_{j}= \frac{1}{(\Delta v)^2}\frac{1}{\sqrt{M_j}}\bigg\{\sqrt{M_j M_{j+1}} \left(\left(\frac{h}{\sqrt{M}}\right)_{j+1}
- \left(\frac{h}{\sqrt{M}}\right)_{j}\right) - \sqrt{M_j M_{j-1}}\left(\left(\frac{h}{\sqrt{M}}\right)_{j} - \left(\frac{h}{\sqrt{M}}\right)_{j-1}\right)\bigg\} \notag\\[4pt]
&\displaystyle\qquad = \frac{1}{(\Delta v)^2}\left(h_{j+1} -\frac{\sqrt{M_{j+1}} + \sqrt{M_{j-1}}}{\sqrt{M_j}}h_j + h_{j-1}\right). 
\end{align*}
Since the new operator $\tilde P$ is symmetric, the Conjugate Gradient (CG) method can be used to get 
$\left(\frac{f}{\sqrt{M}}\right)^{n+1}$. See section 3 in \cite{JinYan} on details for the full discretization. 
\\[2pt]

{\bf Definitions of $\lambda(T)$ and coefficients $D_{ij}^{\varepsilon}$ ($i,j=1,2$)}
\\[6pt]
We recall some definitions given in \cite{Degond2}. In the Boltzmann case, $\lambda(T)$ is given by 
$$\lambda(T) = \frac{2}{3}\int_{\mathbb R^3}\int_{\mathbb S^2} B(v, \Omega)(v, \Omega)^2\, 
M_{0, T}(v)\, d\Omega dv, $$
and in the FPL case, 
$$ \lambda(T) = \frac{2}{3}\int_{\mathbb R^3} B(v)\, M_{0,T}(v)\, dv. $$ 
The coefficients $D_{ij}^{\varepsilon}$ are given by
\begin{align*}
&\displaystyle D_{1j}^{\varepsilon} = \frac{1}{3}\int_{\mathbb R^3} M_{0, T_{\varepsilon}^L}(v)\, \Psi_{j\varepsilon}(|v|)\, |v|^2\, dv, \\[2pt]
&\displaystyle D_{2j}^{\varepsilon} = \frac{1}{6}\int_{\mathbb R^3} M_{0, T_{\varepsilon}^L}(v)\, \Psi_{j\varepsilon}(|v|)\, |v|^4\, dv. 
\end{align*}
$\Psi_i$ is given by the following: The unique solutions $\psi_i^L$, $i=1, 2$, in 
$\left(\text{ker }\Gamma_0^L\right)^{\perp}$, of the equations
$$\Gamma_0^L\, \psi_1^L=v^L, \qquad \Gamma_0^L\, \psi_2^L = \left(\frac{1}{2}\frac{|v^L|^2}{T_0^L}- \frac{3}{2}\right) v^L$$ 
are of the form: $$\psi_i^L = -\Psi_i(|v^L|)\, v^L. $$

\bibliographystyle{siam}
\bibliography{DM.bib}

\begin{thebibliography}{10}

\bibitem{MM-Lemou}
{\sc M.~Bennoune, M.~Lemou, and L.~Mieussens}, {\em Uniformly stable numerical
  schemes for the {B}oltzmann equation preserving the compressible
  {N}avier-{S}tokes asymptotics}, J. Comput. Phys., 227 (2008), pp.~3781--3803.

\bibitem{MC-Bobylev}
{\sc A.~V. Bobylev and I.~F. Potapenko}, {\em Monte {C}arlo methods and their
  analysis for {C}oulomb collisions in multicomponent plasmas}, J. Comput.
  Phys., 246 (2013), pp.~123--144.

\bibitem{BPS}
{\sc A.~V. Bobylev, I.~F. Potapenko, and P.~H. Sakanaka}, {\em Relaxation of
  two-temperature plasma}, Phys. Rev. E, 56 (1997), pp.~2081--2093.

\bibitem{CC}
{\sc C.~Cercignani}, {\em The {B}oltzmann equation and its applications},
  Applied Mathematical Sciences, Springer-Verlag, 67 (1988).

\bibitem{Chapman58}
{\sc S.~Chapman and T.~G. Cowling}, {\em The mathematical theory of non-uniform
  gases}, Cambridge University Press, New York, 1958.

\bibitem{degond2017asymptotic}
{\sc P.~Degond and F.~Deluzet}, {\em Asymptotic-preserving methods and
  multiscale models for plasma physics}, Journal of Computational Physics, 336
  (2017), pp.~429--457.

\bibitem{Degond}
{\sc P.~Degond and B.~Lucquin-Desreux}, {\em The asymptotics of collision
  operators for two species of particles of disparate masses}, Math. Models
  Methods Appl. Sci., 6 (1996), pp.~405--436.

\bibitem{Degond2}
\leavevmode\vrule height 2pt depth -1.6pt width 23pt, {\em Transport
  coefficients of plasmas and disparate mass binary gases}, Transport Theory
  Statist. Phys., 25 (1996), pp.~595--633.

\bibitem{dimarco2011exponential}
{\sc G.~Dimarco and L.~Pareschi}, {\em Exponential {R}unge--{K}utta methods for
  stiff kinetic equations}, SIAM Journal on Numerical Analysis, 49 (2011),
  pp.~2057--2077.

\bibitem{Filbet-Jin}
{\sc F.~Filbet and S.~Jin}, {\em A class of asymptotic-preserving schemes for
  kinetic equations and related problems with stiff sources}, J. Comput. Phys.,
  229 (2010), pp.~7625--7648.

\bibitem{MG18}
{\sc I.~M. Gamba and M.~Pavi\'c-$\breve{C}$oli\'c}, {\em On existence and
  uniqueness to homogeneous {B}oltzmann flows of monatomic gas mixtures},
  preprint (2018).

\bibitem{GHJ}
{\sc C.~Goebel, C.~Harris, and E.~Johnson}, {\em Near normal behavior of
  disparate mass gas mixtures, in {R}arefied gas dynamics}, J. L. Potter (ed.),
  51, Part 1 (1977), pp.~109--122.

\bibitem{Grad}
{\sc H.~Grad}, {\em in {R}arefied {G}as {D}ynamics}, F. Devienne (ed.),
  (1960), pp.~10--138.

\bibitem{Bird}
{\sc J.~Hirschfelder, C.~Curtis, and R.~Bird}, {\em Molecular theory of gases
  and liquids}, Wiley, New-York, 1954.

\bibitem{jin1999efficient}
{\sc S.~Jin}, {\em Efficient asymptotic-preserving ({AP}) schemes for some
  multiscale kinetic equations}, SIAM Journal on Scientific Computing, 21
  (1999), pp.~441--454.

\bibitem{jin2010asymptotic}
\leavevmode\vrule height 2pt depth -1.6pt width 23pt, {\em Asymptotic
  preserving ({AP}) schemes for multiscale kinetic and hyperbolic equations: a
  review}, Lecture notes for summer school on methods and models of kinetic
  theory (M\&MKT), Porto Ercole (Grosseto, Italy),  (2010), pp.~177--216.

\bibitem{JinLi}
{\sc S.~Jin and Q.~Li}, {\em A {BGK}-penalization-based asymptotic-preserving
  scheme for the multispecies {B}oltzmann equation}, Numerical Methods for
  Partial Differential Equations, 29 (2013), pp.~1056--1080.

\bibitem{JPT1}
{\sc S.~Jin, L.~Pareschi, and G.~Toscani}, {\em Diffusive relaxation schemes
  for multiscale discrete-velocity kinetic equations}, SIAM J. Numer. Anal., 35
  (1998), pp.~2405--2439.

\bibitem{JPT2}
\leavevmode\vrule height 2pt depth -1.6pt width 23pt, {\em Uniformly accurate
  diffusive relaxation schemes for multiscale transport equations}, SIAM J.
  Numer. Anal., 38 (2000), pp.~913--936.

\bibitem{jin2010micro}
{\sc S.~Jin and Y.~Shi}, {\em A micro-macro decomposition-based
  asymptotic-preserving scheme for the multispecies boltzmann equation}, SIAM
  Journal on Scientific Computing, 31 (2010), pp.~4580--4606.

\bibitem{JinYan}
{\sc S.~Jin and B.~Yan}, {\em A class of asymptotic-preserving schemes for the
  {F}okker-{P}lanck-{L}andau equation}, J. Comput. Phys., 230 (2011),
  pp.~6420--6437.

\bibitem{Johnson2}
{\sc E.~A. Johnson}, {\em Energy and momentum equations for disparate--mass
  binary gases}, Phys Fluids,  (1973), pp.~45--49.

\bibitem{Johnson1}
\leavevmode\vrule height 2pt depth -1.6pt width 23pt, {\em Description of
  disparate-mass gases}, Archives of Mechanics, 28 (1976), pp.~803--808.

\bibitem{Landau}
{\sc L.D.Landau}, Zh. Eksp. Teor. Fiz,  (1937), p.~203.

\bibitem{Lemou}
{\sc M.~Lemou}, {\em Relaxed micro-macro schemes for kinetic equations}, C. R.
  Math. Acad. Sci. Paris, 348 (2010), pp.~455--460.

\bibitem{li2014exponential}
{\sc Q.~Li and L.~Pareschi}, {\em Exponential runge--kutta for the
  inhomogeneous boltzmann equations with high order of accuracy}, Journal of
  Computational Physics, 259 (2014), pp.~402--420.

\bibitem{liu2017unified}
{\sc C.~Liu and K.~Xu}, {\em A unified gas kinetic scheme for continuum and
  rarefied flows v: multiscale and multi-component plasma transport},
  Communications in Computational Physics, 22 (2017), pp.~1175--1223.

\bibitem{liu2016unified}
{\sc C.~Liu, K.~Xu, Q.~Sun, and Q.~Cai}, {\em A unified gas-kinetic scheme for
  continuum and rarefied flows iv: Full boltzmann and model equations}, Journal
  of Computational Physics, 314 (2016), pp.~305--340.

\bibitem{MHG14}
{\sc A.~Munaf\`o, J.~R. Haack, I.~M. Gamba, and T.~E. Magin}, {\em A
  spectral-{L}agrangian {B}oltzmann solver for a multi-energy level gas}, J.
  Comput. Phys., 264 (2014), pp.~152--176.

\bibitem{Petit75}
{\sc J.~P. Petit and J.~S. Darrozes}, {\em Une nouvelle formulation des
  $\acute{e}$quations du mouvement d'un gaz ionis$\acute{e}$ dans un
  r$\acute{e}$gime domin$\acute{e}$ par les collisions}, Journal de
  M$\acute{e}$canique, 14 (1975), pp.~745--759.

\bibitem{Cheng-Gamba}
{\sc C.~Zhang and I.~Gamba}, {\em A conservative scheme for {V}lasov {P}oisson
  {L}andau modeling collisional plasmas},  (2016).

\end{thebibliography}
\end{document}